\documentclass[12pt,reqno]{amsart}
\usepackage{amsmath,amsfonts,amsthm,amsopn,amssymb}
\usepackage{cite,marginnote}
\usepackage{color,enumitem,graphicx}
\usepackage[colorlinks=true,urlcolor=blue,
citecolor=red,linkcolor=blue,linktocpage,pdfpagelabels,
bookmarksnumbered,bookmarksopen]{hyperref}
\usepackage[english]{babel}

\usepackage[left=2.9cm,right=2.9cm,top=2.8cm,bottom=2.8cm]{geometry}
\usepackage[hyperpageref]{backref}

\numberwithin{equation}{section}

\makeindex
\newtheorem{theorem}{Theorem}[section]
\newtheorem{lemma}{Lemma}[section]
\newtheorem{corollary}{Corollary}[section]
\newtheorem{proposition}{Proposition}[section]
\newtheorem{remark}{Remark}[section]
\newtheorem{example}{Example}[section]
\newtheorem{definition}{Definition}[section]
\newtheorem{theoremletter}{Theorem}

\newcommand{\ud}{\mathrm{d}}
\newcommand{\RN}{\mathbb R^N}
\newcommand{\om}{\Omega}

\newcommand{\iy}{\infty}

\newcommand{\s}{\section}
\newcommand{\dd}{\delta}
\newcommand{\DD}{\Delta}
\newcommand{\g}{\gamma}
\newcommand{\G}{\Gamma}
\newcommand{\na}{\nabla}

\newcommand{\la}{\lambda}
\newcommand{\lb}{\Lambda}

\newcommand{\R}{\mathbb R}
\newcommand{\al}{\alpha}

\newcommand{\ti}{\tilde}

\newcommand{\re}[1]{\eqref{#1}}

\newcommand{\rg}{\rightarrow}

\newcommand{\lan}{\langle}
\newcommand{\ran}{\rangle}
\newcommand{\e}{\varepsilon}
\newcommand{\vp}{\varphi}

\newcommand{\lab}{\label}
\newcommand{\bt}{\begin{theorem}}
\newcommand{\et}{\end{theorem}}
\newcommand{\bl}{\begin{lemma}}
\newcommand{\el}{\end{lemma}}
\newcommand{\bd}{\begin{definition}}
\newcommand{\ed}{\end{definition}}
\newcommand{\bc}{\begin{corollary}}
\newcommand{\ec}{\end{corollary}}
\newcommand{\bp}{\begin{proof}}
\newcommand{\ep}{\end{proof}}
\newcommand{\bx}{\begin{example}}
\newcommand{\ex}{\end{example}}
\newcommand{\bi}{\begin{exercise}}
\newcommand{\ei}{\end{exercise}}
\newcommand{\bo}{\begin{proposition}}
\newcommand{\eo}{\end{proposition}}
\newcommand{\br}{\begin{remark}}
\newcommand{\er}{\end{remark}}
\newcommand{\be}{\begin{equation}}
\newcommand{\ee}{\end{equation}}
\newcommand{\ba}{\begin{align}}
\newcommand{\ea}{\end{align}}
\newcommand{\bn}{\begin{enumerate}}
\newcommand{\en}{\end{enumerate}}
\newcommand{\bg}{\begin{align*}}
\newcommand{\bcs}{\begin{cases}}
\newcommand{\ecs}{\end{cases}}

\newcommand{\bean}{\begin{eqnarray*}}
\newcommand{\eean}{\end{eqnarray*}}

\title[Ground states and semiclassical states]{Ground states and semiclassical states of nonlinear Choquard equations involving Hardy-Littlewood-Sobolev critical growth}

\author[D.\ Cassani]{Daniele Cassani}
\author[J. J.\ Zhang]{Jianjun Zhang}

\address[D.\ Cassani]{\newline\indent Dip. di Scienza e Alta Tecnologia
\newline\indent
Universit\`{a} degli Studi dell'Insubria
\newline\indent
via Valleggio 11, 22100 Como,Italy
\newline\indent
and 
\newline\indent
RISM - Riemann International School of Mathematics
\newline\indent
via G.B. Vico 46, 21100 Varese, Italy}
\email{\href{mailto:Daniele.Cassani@uninsubria.it}{Daniele.Cassani@uninsubria.it}}

\address[J. J.\ Zhang]{\newline\indent College of Mathematics and Statistics
\newline\indent
Chongqing Jiaotong University
\newline\indent
Chongqing 400074, PR China
\newline\indent and
\newline\indent Dip. di Scienza e Alta Tecnologia
\newline\indent
Universit\`{a} degli Studi dell'Insubria
\newline\indent
via Valleggio 11, 22100 Como,Italy}
\email{\href{mailto:zhangjianjun09@tsinghua.org.cn}{zhangjianjun09@tsinghua.org.cn}}

\thanks{J.J. Zhang was partially supported by the Science Foundation of Chongqing Jiaotong University(15JDKJC-B033).}

\subjclass[2000]{35B25, 35B33, 35J61}
\keywords{Ground states, Semiclassical states, Nonlocal equations, Hardy-Littlewood-Sobolev inequality}

\begin{document}

\begin{abstract}
We are concerned with the existence of ground states for nonlinear Choquard equations involving a critical nonlinearity in the sense of Hardy-Littlewood-Sobolev. Our result complements previous results by Moroz and Van Schaftingen where the subcritical case was considered. Then, we focus on the existence of semi-classical states and by using a truncation argument approach, we establish the existence and concentration of single peak solutions concentrating around minima of the Schr\"odinger potential, as the Planck constant goes to zero. The result is robust in the sense that the nonlinearity is not required to satisfy {\it monotonicity} conditions nor the {\it Ambrosetti-Rabinowitz} condition.
\end{abstract}

\maketitle

\s{Introduction}
\renewcommand{\theequation}{1.\arabic{equation}}

\noindent This paper is concerned with the following nonlinear Choquard equation
\begin{equation}\label{q1}
-\e^2\DD v+V(x)v=\e^{-\al}(I_\al\ast F(v))f(v),\,\, v>0,\,\,x\in\RN,
\end{equation}
where $N\geq 3$, $\alpha\in(0,N)$, $F$ is the primitive function of $f$, $I_\al$ is the Riesz potential defined for every $x\in\RN\setminus\{0\}$ by
$$
I_\al(x):=\frac{A_\al}{|x|^{N-\al}},\,\,\mbox{where}\,\, A_\al=\frac{\G((N-\al)/2)}{\G(\al/2)\pi^{N/2}2^\al},\hbox{ and $\G$ is the Gamma function.}
$$
The Schr\"odinger potential $V$ satisfies the following 
\begin{itemize}

\item [(V1)] $V\in C(\RN,\R)$ and $\inf_{x\in\RN}V(x)>0$.

\end{itemize}
When $\e=1$ and $V(x)=a>0$,  \re{q1} reduces to the following nonlocal elliptic equation
\begin{eqnarray}\label{lb1}
-\Delta u+au=(I_\al\ast F(u))f(u), \ x\in \RN,
\end{eqnarray}
which is variational, in the sense that, $u$ is a solution to \re{lb1} if and only if $u$ is a critical point of the following energy functional
$$
L_a(u)=\frac{1}{2}\int_{\RN}|\na u|^2+au^2-\frac{1}{2}\int_{\RN} (I_\al\ast F(u))F(u),\ \ u\in H^1(\RN).
$$
In the relevant physical case in which $N=3$, $\al=2$ and $F(s)=s^2/2$, \re{lb1} turns into the equation
\begin{eqnarray}\label{limit equation2}
-\Delta u+au=(I_2\ast u^2)u, \ x\in \R^3,
\end{eqnarray}
which is the so-called Choquard-Pekar equation and can goes back to the description of a polaron at rest in Quantum Field Theory by S.I. Pekar \cite{Pekar}. Moreover, if $u$ is a solution of \re{limit equation2}, then $\psi(x,t)=e^{it}u(x)$ is a solitary wave of the focusing time-dependent Hartree equation
$$
i\psi_t=-\DD\psi-(I_2\ast \psi^2)\psi, \ (t,x)\in\R^+\times\R^3.
$$
In 1976, P. Choquard introduced this type of equations to describe an electron trapped in its own hole as an approximation to Hartree-Fock theory for a one component plasma, see \cite{LS}. It also arises in multiple particles systems\cite{Gross, LS} and quantum mechanics\cite{Penrose1,Penrose2,Penrose3}. In the pioneering work \cite{Lieb1}, E.H. Lieb first proved the existence and uniqueness of positive solutions to \re{limit equation2}. Later, multiplicity results for \eqref{limit equation2} were obtained by P.L. Lions\cite{Lions2,Lions3} by variational methods.

For $F(u)=|u|^p/p$, \re{lb1} can be reduced to the following general stationary Hartree equation
\begin{eqnarray}\label{limit equation3}
-\Delta u+u=(I_\al\ast |u|^{p}) |u|^{p-2}u, \ x\in\RN.
\end{eqnarray}
Let us recall that in the local case
\begin{eqnarray}\label{limit equation333}
-\Delta u+u=|u|^{p-2}u, \ x\in\RN,
\end{eqnarray}
it is well known that positive solutions with finite energy are radially symmetric, unique and non-degenerate, see \cite{Gidas,Oh}. Here, in contrast to the local problem \re{limit equation333}, the standard moving planes approach seems to be unsettled at the moment to deal with the nonlocal version \re{limit equation3}. The classification of positive solutions to \re{limit equation3}, even in the particular case $p=2$, has been remained a longstanding open problem. By using an integral version of the moving planes method, introduced by W. Chen et al. \cite{Chen}, L. Ma and L. Zhao \cite{MZ} gave a breakthrough on this problem. With some restrictions on $\alpha$, $p$ and $N$, they proved that positive solutions to \eqref{limit equation3} are, up to translations, radially symmetric and unique. In \cite{MV3}, V. Moroz and J. van Schaftingen further improved the result in of \cite{MZ} by establishing the existence of ground state solutions to \eqref{limit equation3} within an optimal range of $p$. More recently, V. Moroz and J. van Schaftingen\cite{MV1} considered the more general Choquard equation \eqref{lb1} and, in the spirit of Berestycki and Lions, obtained the existence of ground state solutions with sufficient and almost necessary conditions on the nonlinearity $f$. For more details on this subject, we refer to the survey \cite{MV4}.

In the above quoted literature, only the subcritical case was considered. The first purpose of the present work is to investigate the existence of ground state solutions to \re{lb1} involving critical growth in the sense of Hardy-Littlewood-Sobolev inequality.

\bd
$u$ is said to be a ground state solution of \re{lb1} if $u$ is a solution of \re{lb1} with the least energy $L_a$ among all nontrivial solutions to \re{lb1}.
\ed
\noindent Throughout this paper we assume $f\in C(\R^+,\R)$ which satisfies
\begin{itemize}
\item [(F1)] $\lim_{t\rg 0^+}f(t)/t=0$;
\item [(F2)] $\lim_{t\rg +\infty}f(t)/t^{\frac{\al+2}{N-2}}=1$;
\item [(F3)] there exist $\mu>0$ and $q\in(2,(N+\al)/(N-2))$ such that $$f(t)\ge t^{(2+\al)/(N-2)}+\mu t^{q-1},\,\,t>0.$$
\end{itemize}
Our first main result is the following
\bt\lab{Theorem 1} Let $\al\in((N-4)_+,N)$ and 
$$
q>\max\left\{1+\frac{\al}{N-2},\frac{N+\al}{2(N-2)}\right\},
$$
and assume conditions $(F1)$--$(F3)$. Then, for any $a>0$, \re{lb1} admits a ground state solution.
\et
\br Let us point out that to ensure the existence of ground states to \re{lb1}, the assumption $(F3)$ plays a crucial role.
Without $(F3)$, the assumptions $(F1)$-$(F2)$ can not guarantee the existence of ground states to \re{lb1}. Here we give a counterexample: let $\al\in((N-4)_+,N)$ and $f(t)=|t|^{(4+\al-N)/(N-2)}t$, which satisfies $(F1)$-$(F2)$ but not $(F3)$. By a Pohoz\v{a}ev's type identity (see Lemma \ref{Pohozaev}, Section 2),  \re{lb1} has no nontrivial solutions.
\er

The second purpose of this paper is to investigate the profile of positive solutions to \re{q1} as the adimensionalized Planck constant $\e\rg0$, whose motivation goes back to the pioneering work of A. Floer and A. Weinstein \cite{F-W} (see also \cite{Oh}) concerning the Schr\"odinger equation
\be\lab{sch}
-\e^2\DD u+V(x)u=f(u),\,\,\ x\in\RN.
\ee
An interesting class of solutions to \re{sch} are families of solutions which develop a spike shape around some point in $\RN$ as $\e\rg 0$. From the physical point of view, for $\e>0$ small, these solutions give the so-called semi-classical states, which describe the transition from quantum mechanics to classical mechanics. For the detailed physical background, we refer to \cite{Oh} and references therein, see also \cite{evans_zworski}. By a Lyapunov-Schmidt reduction approach based on a non-degeneracy condition, in \cite{F-W,Oh}, the authors obtained the existence of solutions to \re{sch} exhibiting a single peak or multi peaks concentrating, as $\e\rg0$, around any given non-degenerate critical points of the potential $V$. However, the non-degeneracy condition holds only for a restricted class of nonlinearities $f$. In the last decade, considerable attention has been paid to relax or remove the non-degeneracy condition in the singularly perturbed problems. By using a variational approach, P.H. Rabinowitz \cite{Rab} obtained the existence of positive solutions to \re{sch} for small $\e>0$ by assuming the following global potential well condition
$$
\liminf_{|x|\rg\iy}V(x)>\inf_{\RN}V(x).
$$
Later, by using a penalization approach, M. del Pino and P. Felmer \cite{DF} weakened the global potential well condition above to the the following local condition
\begin{itemize}

\item [(V2)] there is a bounded domain $O\subset\RN$ such that $$0<m\equiv\inf_{x\in O}V(x)<\min_{x\in \partial O}V(x),$$

\end{itemize}
and proved the existence of a single-peak solution to \re{sch}. In \cite{Rab, DF}, the non-degeneracy condition is not required. Some related results can be found in \cite{WX,Felmer2,Pino,DPR,Alves1} and the references therein.
In 2007, J. Byeon and L. Jeanjean \cite{byeon} introduced a new penalization approach and constructed a spike solution under hypothesis $(V2)$ and almost optimal hypotheses on $f$: namely, the Berestycki-Lions\ conditions \cite{Lions}. For further references, we refer the reader to \cite{BT1,BT2,byeon4} for the subcritical case, and to \cite{zhang-chen-zou, ZZ} for the critical case. 

\vskip0.1in
\noindent We state the second main result of this paper as follows.

\bt\lab{Theorem 2} Assume $(V1)$-$(V2)$ and the assumptions of Theorem \ref{Theorem 1}. Let $\mathcal{M}\equiv\{x\in O: V(x)=m\}$. Then, for small $\e>0$, \re{q1} admits a positive solution $v_{\e}$, which satisfies:
\begin{itemize}
\item [(i)] there exists a local maximum point $x_\e\in O$ of $v_\e$ such that $$\lim_{\e\rg 0}dist(x_\e,\mathcal{M})=0,$$ and $w_\e(x)\equiv v_\e(\e x+x_\e)$ converges (up to a subsequence) uniformly to a ground state solution of
$$ -\DD u+mu=(I_\al\ast F(u))f(u),\ u>0,\ u\in H^1(\RN);$$
\item [(ii)] $v_\e(x)\le C\exp(-\frac{c}{\e}|x-x_\e|)$ for some $c,C>0$.
\end{itemize}
\et
In \cite{WW}, J. Wei and M. Winter considered the Schr\"{o}dinger-Newton system
\begin{equation}\label{main equation5}
-\e^2\DD v+V(x)v=\e^{-2}(I_2\ast v^2)v,\ \ x\in \R^3
\end{equation}
and by using a Lyapunov-Schmidt reduction method and with the assumption $(V1)$, proved the existence of multi-bump solutions concentrating around local minima, local maxima or non-degenerate critical points of $V$. Let us mention that when $(V1)$ fails to hold and the potential vanishes somewhere, the problem becomes more difficult. In \cite{Secchi}, S. Secchi considered the Schr\"{o}dinger-Newton system \re{main equation5} with a positive decaying electric potential and by virtue of perturbative methods, proved the existence and concentration of bound states near local minima (or maxima) points of $V$ as $\e\rg0$. Recently, by a nonlocal penalization technique, V. Moroz and J. Van Schaftingen \cite{MV2} obtained a family of single spike solutions of the Choquard equation
$$
-\e^2\DD v+V(x)v=\e^{-2}(I_\al\ast |v|^p)|v|^{p-2}v,\ \ x\in \RN.
$$
around the local minimum of $V$ as $\e\rg0$. Moreover, in \cite{MV2}, the assumption on the decay of $V$ and the admissible range for $p\ge2$ are optimal. More recently, adopting the penalization argument introduced in \cite{byeon}, M. Yang et al. \cite{YZZ} investigated the existence and concentration of solutions to \re{q1} under the local potential well condition $(V2)$ and a mild assumption on $f$. In particular, the {\it Ambrosetti-Rabinowtiz} condition and the {\it monotonicity} condition on $f(t)/t$ are not required. For more related results, we refer to \cite{alves,Sun,CSS, MN, Nolasco, Secchi,YD,Squa} and the references therein. However, the above quoted result cover the subcritical case and the critical case, in the terms of the Hardy-Littlewood-Sobolev inequality, remained open. In \cite{alves0}, C. O. Alves et al. considered the ground state solutions of the Choquard equation \re{q1} in $\R^2$. By variational methods, the authors proved the existence and concentration of ground states to \re{q1} involving critical exponential growth in the sense of Trudinger-Moser. A natural open problem which has not been settled before is whether \re{q1} develop similar concentration phenomena in the case of critical growth. Here we answer this question completing the study carried out in the above quoted literature. 
\vskip0.1in
\noindent In Section 2, we use a monotonicity trick due to Jeanjean \cite{Jeanjean} and a suitable decomposition of Palais-Smale sequences to prove Theorem \ref{Theorem 1}. In Section 3, we prove Theorem \ref{Theorem 2} by means of a truncation approach.


\s{Proof of Theorem \ref{Theorem 1}}

\renewcommand{\theequation}{2.\arabic{equation}}

\noindent
In this section, we are concerned with the existence of ground state solutions to \re{lb1}. Let $a>0$ and denote the least energy of \re{lb1} by
$$
E_a=\inf\left\{L_a(u): L_a'(u)=0\,\,\hbox{in $H^{-1}(\RN)$},\, \, u\in H^1(\RN)\setminus\{0\}\right\}.
$$
In what follows, let $H^1(\RN)$ be endowed with the norm
$$
\|u\|=\left(\int_{\RN}|\na u|^2+a|u|^2\right)^2,\,\, u\in H^1(\RN).
$$

\noindent Before proving Theorem \ref{Theorem 1}, we introduce some preliminary results. First, the following Hardy-Littlewood-Sobolev inequality will be used frequently later.
\bl\lab{hls}{\rm \cite[Theorem 4.3]{LL}}
Let $s, r>1$ and $0<\al<N$ with $1/s+1/r=1+\al/N$, $f\in L^s(\RN)$ and $g\in L^r(\RN)$, then there exists a positive constant $C(s, N, \al)$ (independent of $f, g$) such that
$$
\left|\int_{\RN}\int_{\RN}f(x)|x-y|^{\al-N}g(y)\,\ud x\ud y\right|\le C(s, N, \al)\|f\|_s\|g\|_r.
$$
In particular, if $s=r=2N/(N+\al)$, the sharp constant
$$
\mathcal{C}_\al:=\pi^{\frac{N-\al}{2}}\frac{\G(\al/2)}{\G((N+\al)/2)}\left[\frac{\G(N/2)}{\G(N)}\right]^{-\al/N}.
$$
\el
\br\lab{rr1}
By the Hardy-Littlewood-Sobolev inequality above, for any $u\in L^s(\RN)$ with $s\in(1,N/\al)$, $I_\al\ast v\in L^{Ns/(N-\al s)}(\RN)$. Moreover, $I_\al\in\mathcal{L}(L^s(\RN),L^{Ns/(N-\al s)}(\RN))$ and
$$
\|I_\al\ast v\|_{\frac{Ns}{N-\al s}}\le C(s, N,\al)\|v\|_s.
$$
\er

\subsection{Brezis-Lieb lemma and Splitting lemma}
In this section, we give a Brezis-Lieb lemma and splitting lemma for the nonlocal term of the functional.
\bl\lab{Lieb} {\rm (Brezis-Lieb Lemma)} Assume $\al\in(0,N)$ and there exists $C>0$ such that
$$
|f(t)|\le C(|t|^{\frac{\al}{N}}+|s|^{\frac{\al+2}{N-2}}),\,\, s\in\R.
$$
Let $\{u_n\}\subset H^1(\RN)$ such that $u_n\rg u$ weakly in $H^1(\RN)$ and a.e. in $\RN$ as $n\rg\iy$, then
$$
\int_{\RN}(I_\al\ast F(u_n))F(u_n)=\int_{\RN}(I_\al\ast F(u_n-u))F(u_n-u)+\int_{\RN}(I_\al\ast F(u))F(u)+o_n(1),
$$
where $o_n(1)\rg0$ as $n\rg\iy$.
\el
To prove Lemma \ref{Lieb}, we recall the following lemma, which states that pointwise convergence of a bounded sequence implies weak convergence.
\bl\lab{almost}{\rm\cite[ Theorem 4.2.7]{Willem}}
Let $\om\subset\RN$ be a domain and $\{u_n\}$ be bounded in $L^q(\om)$ with some $q>1$. Then if $u_n\rg u$ a.e. in $\om$ as $n\rg\iy$, then $u_n\rg u$ weakly in $L^q(\om)$ as $n\rg\iy$
\el
\noindent{\it Proof of Lemma \ref{Lieb}.}
Observe that
\begin{align*}
&\int_{\RN}(I_\al\ast F(u_n))F(u_n)-(I_\al\ast F(u_n-u))F(u_n-u)-(I_\al\ast F(u))F(u)\\
&=\int_{\RN}(I_\al\ast[F(u_n)+F(u_n-u)])[F(u_n)-F(u_n-u)]-(I_\al\ast F(u))F(u),
\end{align*}
and there exists $C>0$ such that $|F(s)|\le C(|s|^{(N+\al)/N}+|s|^{(N+\al)/(N-2)})$ for all $s\in\R$, which implies $F(u)\in L^{2N/(N+\al)}(\RN)$. For any $\dd>0$ small, by the Hardy-Littlewood-Sobolev inequality, there exists $K_1>0$ such that
$$
\left|\int_{\om_1}(I_\al\ast F(u))F(u)\right|\le\dd/6,\,\,\om_1:=\{x\in\RN: |u(x)|\ge K_1\}.
$$
Meanwhile, by the Hardy-Littlewood-Sobolev inequality,
{\allowdisplaybreaks
\begin{align*}
&\left|\int_{\om_1}(I_\al\ast[F(u_n)+F(u_n-u)])[F(u_n)-F(u_n-u)]\right|\\
&\le C\left(\int_{\RN}|F(u_n)+F(u_n-u)|^{\frac{2N}{N+\al}}\right)^{\frac{N+\al}{2N}}
\left(\int_{\om_1}|F(u_n)-F(u_n-u)|^{\frac{2N}{N+\al}}\right)^{\frac{N+\al}{2N}}\\
&\le C(N,\al)\left(\int_{\om_1}|F(u_n)-F(u_n-u)|^{\frac{2N}{N+\al}}\right)^{\frac{N+\al}{2N}}.
\end{align*}}%
Here we used the fact that $\{u_n\}$ is bounded in $H^1(\RN)$. It is easy to see there exists $c>0$ such that
\begin{align*}
&|F(u_n)-F(u_n-u)|^{\frac{2N}{N+\al}}\\
&\le c(|u_n|^{\frac{2\al}{N+\al}}|u|^{\frac{2N}{N+\al}}+|u_n|^{\frac{2+\al}{N-2}\frac{2N}{N+\al}}|u|^{\frac{2N}{N+\al}}+u^2+|u|^{\frac{2N}{N-2}}),\,\, x\in\RN.
\end{align*}
Then by H\"older's inequality,
$$
\int_{\om_1}|u_n|^{\frac{2\al}{N+\al}}|u|^{\frac{2N}{N+\al}}
\le\left(\int_{\om_1}u_n^2\right)^{\frac{\al}{N+\al}}\left(\int_{\om_1}u^2\right)^{\frac{N}{N+\al}}
$$
and
$$
\int_{\om_1}|u_n|^{\frac{2+\al}{N-2}\frac{2N}{N+\al}}|u|^{\frac{2N}{N+\al}}
\le\left(\int_{\om_1}|u_n|^{\frac{2N}{N-2}}\right)^{\frac{2+\al}{N+\al}}\left(\int_{\om_1}|u|^{\frac{2N}{N-2}}\right)^{\frac{N-2}{N+\al}}.
$$
So for $\dd$ given above and $K_1$ fixed but large enough, we get for any $n$,
$$
\left|\int_{\om_1}(I_\al\ast[F(u_n)+F(u_n-u)])[F(u_n)-F(u_n-u)]\right|\le \dd/6.
$$
Similarly, let $\om_2:=\{x\in\RN: |x|\ge R\}\setminus\om_1$ with $R>0$ large enough, we have for any $n$,
$$
\left|\int_{\om_2}(I_\al\ast F(u))F(u)\right|\le\dd/6
$$
and
$$
\left|\int_{\om_2}(I_\al\ast[F(u_n)+F(u_n-u)])[F(u_n)-F(u_n-u)]\right|\le \dd/6.
$$
For $K_2>K_1$, let $\om_3(n):=\{x\in\RN: |u_n(x)|\ge K_2\}\setminus(\om_1\cup\om_2)$, then if $\om_3(n)\not=\emptyset$, we know $|u(x)|<K_1$ and $|x|<R$ for any $x\in\om_3(n)$. Noting that $u_n\rg u$ a.e. in $\om$ as $n\rg\iy$, then it follows from the Severini-Egoroff theorem that $u_n$ converges to $u$ in measure in $B_R(0)$, which implies that $|\om_3(n)|\rg0$ as $n\rg\iy$. Then similar as above, we have for $n$ large enough,
$$
\left|\int_{\om_3(n)}(I_\al\ast F(u))F(u)\right|\le\dd/6
$$
and
$$
\left|\int_{\om_3(n)}(I_\al\ast[F(u_n)+F(u_n-u)])[F(u_n)-F(u_n-u)]\right|\le \dd/6.
$$
Finally, we estimate the terms
$$
\int_{\om_4(n)}(I_\al\ast[F(u_n)+F(u_n-u)])[F(u_n)-F(u_n-u)]-(I_\al\ast F(u))F(u),
$$
where $\om_4(n)=\RN\setminus(\om_1\cup\om_2\cup\om_3(n))$. Obviously, $\om_4(n)\subset B_R(0)$. By Lebesgue's dominated convergence theorem, $$\lim_{n\rg\iy}\int_{\om_4(n)}|F(u_n-u)|^{\frac{2N}{N+\al}}=0,\,\, \lim_{n\rg\iy}\int_{\om_4(n)}|F(u_n)-F(u)|^{\frac{2N}{N+\al}}=0,$$
which implies by the Hardy-Littlewood-Sobolev inequality that as $n\rg\iy$,
\begin{align*}
&\left|\int_{\om_4(n)}(I_\al\ast[F(u_n)+F(u_n-u)])F(u_n-u)\right|\\
&\le C(N,\al)\left(\int_{\om_4(n)}|F(u_n-u)|^{\frac{2N}{N+\al}}\right)^{\frac{N+\al}{2N}}\rg0
\end{align*}
and
\begin{align*}
&\left|\int_{\om_4(n)}(I_\al\ast[F(u_n)+F(u_n-u)])[F(u_n)-F(u)]\right|\\
&\le C(N,\al)\left(\int_{\om_4(n)}|F(u_n)-F(u)|^{\frac{2N}{N+\al}}\right)^{\frac{N+\al}{2N}}\rg0
\end{align*}
Then let $H_n=F(u_n)+F(u_n-u)-F(u)$, we have
\begin{align*}
&\lim_{n\rg\iy}\int_{\om_4(n)}(I_\al\ast[F(u_n)+F(u_n-u)])[F(u_n)-F(u_n-u)]-(I_\al\ast F(u))F(u)\\
&=\lim_{n\rg\iy}\int_{\om_4(n)}(I_\al\ast H_n)F(u).
\end{align*}
Noting that ${H_n}$ is bounded in $L^{2N/(N+\al)}(\RN)$ and $H_n\rg0$ a. e. in $\RN$ as $n\rg\iy$, by Lemma \ref{almost}, $H_n\rg0$ weakly in $L^{2N/(N+\al)}(\RN)$ as $n\rg\iy$. By Remark \ref{rr1}, $I_\al\ast H_n\rg0$ weakly in $L^{2N/(N-\al)}(\RN)$ as $n\rg\iy$, which yields that
$$
\lim_{n\rg\iy}\int_{\om_4(n)}(I_\al\ast H_n)F(u)=0.
$$
Thus,
$$
\limsup_{n\rg\iy}\left|\int_{\RN}(I_\al\ast F(u_n))F(u_n)-(I_\al\ast F(u_n-u))F(u_n-u)-(I_\al\ast F(u))F(u)\right|\le\dd.
$$
By the arbitrary choice of $\dd$, the proof is completed.
\qed
\vskip0.1in
Next, we give the following splitting lemma.

\bl\lab{Split} {\rm (Splitting Lemma)} Assume $\al\in((N-4)_+,N)$, $(F1)$-$(F2)$ and let $\{u_n\}\subset H^1(\RN)$ such that $u_n\rg u$ weakly in $H^1(\RN)$ and a.e. in $\RN$ as $n\rg\iy$, then passing to a subsequence, if necessary,
$$
\int_{\RN}\big([I_\al\ast F(u_n)]f(u_n)-[I_\al\ast F(u_n-u)]f(u_n-u)-[I_\al\ast F(u)]f(u)\big)\phi=o_n(1)\|\phi\|,
$$
where $o_n(1)\rg 0$ uniformly for any $\phi\in C_0^\iy(\RN)$ as $n\rg\iy$.
\el
To prove Lemma \ref{Split}, we give Lemma \ref{Brezis} and Lemma \ref{Asymptotic} as follows.

\bl\lab{Brezis}
Let $\om\subset\RN$ be a domain and $\{u_n\}\subset H^1(\om)$ such that $u_n\rg u$ weakly in $H^1(\om)$ and a.e. in $\om$ as $n\rg\iy$.
\begin{itemize}
\item [$(i)$] For any $1<q\le r\le 2N/(N-2)$ and $r>2$,
$$
\lim_{n\rg\iy}\int_{\om}\big||u_n|^{q-1}u_n-|u_n-u|^{q-1}(u_n-u)-|u|^{q-1}u\big|^{\frac{r}{q}}=0.
$$
\item [$(ii)$] Assume $h\in C(\R,\R)$ and $h(t)=o(t)$ as $t\rg0$, $|h(t)|\le c(1+|t|^{q})$ for any $t\in\R$ where $q\in(1,(N+2)/(N-2)]$, then
    \begin{itemize}
      \item [$(ii)_1$] for any $r\in [q+1,2N/(N-2)]$,
       $$
         \lim_{n\rg\iy}\int_{\om}\big|H(u_n)-H(u_n-u)-H(u)\big|^{\frac{r}{q+1}}=0,
       $$
where $H(t)=\int_0^th(s)\,\ud s$,
      \item [$(ii)_2$] if we further assume that $\om=\RN$, $\al\in((N-4)_+,N)$ and $\lim_{|t|\rg\iy}h(t)|t|^{-\frac{\al+2}{N-2}}=0$, then
      $$
      \int_{\RN}\left|h(u_n)-h(u_n-u)-h(u)\right|^{\frac{2N}{N+\al}}|\phi|^{\frac{2N}{N+\al}}=o_n(1)\|\phi\|^{\frac{2N}{N+\al}},
      $$
      where $o_n(1)\rg 0$ uniformly for any $\phi\in C_0^\iy(\RN)$ as $n\rg\iy$.
    \end{itemize}
\end{itemize}
\el
\bp
The proofs of $(i)$ and $(ii)_1$ are similar to \cite[Lemma 2.5]{JZ}. We only give the proof of $(ii)_2$ which is inspired by \cite{Weth2} and \cite[Lemma 4.7]{ZZ}.

In the following, let $C$ be positive constants (independent of $\e,k$), which may change from line to line. For any fixed $\e\in(0,1)$, there exists $s_0=s_0(\e)\in (0,1)$ such that
$|h(t)|\le \e |t|$ for $|t|\le 2s_0$. Choosing
$s_1=s_1(\e)>2$ such that $|h(t)|\le \e|t|^{(2+\al)/(N-2)}$ for
$|t|\ge s_1-1$. From the continuity of $h$, there exists
$\dd=\dd(\e)\in(0,s_0)$ such that $|h(t_1)-h(t_2)|\le s_0\e$ for
$|t_1-t_2|\le \dd, |t_1|,|t_2|\le s_1+1$. Moreover, there exists $c(\e)>0$ such that $|h(t)|\le
c(\e)|t|+\e|t|^{(2+\al)/(N-2)}$ for $t\in \R$. Noting that $\al\in((N-4)_+,N)$, we know $2<4N/(N+\al)<2N/(N-2)$. Then there exists
$R=R(\e)>0$ such that
{\allowdisplaybreaks
\begin{align}\lab{bianhao1}
\int_{\RN\setminus B(0,R)}|h(u)\phi|^{\frac{2N}{N+\al}}&\le C\int_{\RN\setminus B(0,R)}\left(|u|^\frac{2N}{N+\al}+\e|u|^{\frac{2+\al}{N-2}\frac{2N}{N+\al}}\right)|\phi|^{\frac{2N}{N+\al}}\nonumber\\
&\le C\left(\int_{\RN\setminus B(0,R)}|u|^\frac{4N}{N+\al}\right)^{\frac{1}{2}}\left(\int_{\RN}|\phi|^\frac{4N}{N+\al}\right)^{\frac{1}{2}}\\
&\ \ \ +C\e\left(\int_{\RN\setminus B(0,R)}|u|^\frac{2N}{N-2}\right)^{\frac{2+\al}{N+\al}}\left(\int_{\RN}|\phi|^{\frac{2N}{N-2}}\right)^{\frac{N-2}{N+\al}}\nonumber\\
&\le C\e\|\phi\|^{\frac{2N}{N+\al}}.\nonumber
\end{align}}%
Setting $A_n:=\{x\in\RN\setminus B(0,R): |u_n(x)|\le
s_0\}$, then
\begin{align*}
&\int_{A_n\cap\{|u|\le\dd\}}|h(u_n)-h(u_n-u)|^{\frac{2N}{N+\al}}|\phi|^{\frac{2N}{N+\al}}\\
&\le C\e\int_{\RN}\left(|u_n|^{\frac{2N}{N+\al}}+|u_n-u|^{\frac{2N}{N+\al}}\right)|\phi|^{\frac{2N}{N+\al}}\\
&\le C\e\|\phi\|^{\frac{2N}{N+\al}}.
\end{align*}
Let
$B_n:=\{x\in\RN\setminus B(0,R): |u_n(x)|\ge s_1\}$, then
\begin{align*}
&\int_{B_n\cap\{|u|\le\dd\}}|h(u_n)-h(u_n-u)|^{\frac{2N}{N+\al}}|\phi|^{\frac{2N}{N+\al}}\\
&\le C\e\int_{\RN}\left(|u_n|^{\frac{2+\al}{N-2}\frac{2N}{N+\al}}+|u_n-u|^{\frac{2+\al}{N-2}\frac{2N}{N+\al}}\right)|\phi|^{\frac{2N}{N+\al}}\\
&\le C\e\|\phi\|^{\frac{2N}{N+\al}}.
\end{align*}
Setting $C_n:=\{x\in\RN\setminus B(0,R): s_0\le|u_n(x)|\le s_1\}$, then $|C_n|<\iy$ and
{\allowdisplaybreaks
\begin{align*}
&\int_{C_n\cap\{|u|\le\dd\}}|h(u_n)-h(u_n-u)|^{\frac{2N}{N+\al}}|\phi|^{\frac{2N}{N+\al}}\\
&\le(s_0\e)^{\frac{2N}{N+\al}}\int_{C_n\cap\{|u|\le\dd\}}|\phi|^{\frac{2N}{N+\al}}\le(s_0\e)^{\frac{2N}{N+\al}}|C_n|^{\frac{1}{2}}\left(\int_{\RN}|\phi|^{\frac{4N}{N+\al}}\right)^{\frac{1}{2}}\\
&\le\e^{\frac{2N}{N+\al}}\left(\int_{C_n}|u_n|^{\frac{4N}{N+\al}}\right)^{\frac{1}{2}}\left(\int_{\RN}|\phi|^{\frac{4N}{N+\al}}\right)^{\frac{1}{2}}\le C\e\|\phi\|^{\frac{2N}{N+\al}}.
\end{align*}}%
Thus, $(\RN\setminus B(0,R))\cap\{|u|\le\dd\}=A_n\cup B_n\cup C_n$ and
$$ \int_{(\RN\setminus B(0,R))\cap\{|u|\le\dd\}}|h(u_n)-h(u_n-u)|^{\frac{2N}{N+\al}}|\phi|^{\frac{2N}{N+\al}}\le
 C\e\|\phi\|^{\frac{2N}{N+\al}}\ \ \mbox{for all}\ \ n.
$$
Obviously, for $\e$ given above, there exists $c(\e)>0$ such that
\begin{align*}
|h(u_n)-h(u_n-u)|^{\frac{2N}{N+\al}}&\le
\e(|u_n|^{\frac{2+\al}{N-2}\frac{2N}{N+\al}}+|u_n-u|^{\frac{2+\al}{N-2}\frac{2N}{N+\al}})\\
& \ \ \ +c(\e)(|u_n|^{\frac{2N}{N+\al}}+|u_n-u|^{\frac{2N}{N+\al}})
\end{align*}
and
{\allowdisplaybreaks
\begin{align*}
&\int_{(\RN\setminus B(0,R))\cap\{|u|\ge\dd\}}|h(u_n)-h(u_n-u)|^{\frac{2N}{N+\al}}|\phi|^{\frac{2N}{N+\al}}\\
& \le \int_{(\RN\setminus
B(0,R))\cap\{|u|\ge\dd\}}\e(|u_n|^{\frac{2+\al}{N-2}\frac{2N}{N+\al}}+|u_n-u|^{\frac{2+\al}{N-2}\frac{2N}{N+\al}})|\phi|^{\frac{2N}{N+\al}}\\
&\quad +c(\e)(|u_n|^{\frac{2N}{N+\al}}+|u_n-u|^{\frac{2N}{N+\al}})|\phi|^{\frac{2N}{N+\al}}\\
& \le C\e\|\phi\|^{\frac{2N}{N+\al}}+c(\e)\int_{(\RN\setminus
B(0,R))\cap\{|u|\ge\dd\}}(|u_n|^{\frac{2N}{N+\al}}+|u_n-u|^{\frac{2N}{N+\al}})|\phi|^{\frac{2N}{N+\al}}.
\end{align*}}%
Noting that $0<\al+4-N<N+\al$ and $|(\RN\setminus B(0,R))\cap\{|u|\ge\dd\}|\rg 0$ as $R\rg\iy$, there exists $R=R(\e)$ large enough, such that
\begin{align*}
&\int_{(\RN\setminus
B(0,R))\cap\{|u|\ge\dd\}}c(\e)(|u_n|^{\frac{2N}{N+\al}}+|u_n-u|^{\frac{2N}{N+\al}})|\phi|^{\frac{2N}{N+\al}}\\
&\le c(\e)\left[\left(\int_{\RN}|u_n|^{\frac{2N}{N-2}}\right)^{\frac{N-2}{N+\al}}+\left(\int_{\RN}|u_n-u|^{\frac{2N}{N-2}}\right)^{\frac{N-2}{N+\al}}\right]\\
&\ \ \ \ \times\left(\int_{\RN}|\phi|^{\frac{2N}{N-2}}\right)^{\frac{N-2}{N+\al}}|(\RN\setminus
B(0,R))\cap\{|u|\ge\dd\}|^{\frac{\al+4-N}{N+\al}}\\
& \le  \e\|\phi\|^{\frac{2N}{N+\al}}.
\end{align*}
Then for any $n$,
$$
\int_{(\RN\setminus
B(0,R))\cap\{|u|\ge\dd\}}|h(u_n)-h(u_n-u)|^{\frac{2N}{N+\al}}|\phi|^{\frac{2N}{N+\al}}\le C\e\|\phi\|^{\frac{2N}{N+\al}}.
$$
Thus, by \re{bianhao1}, for any $n$,
\be\lab{bianhao5}
\int_{\RN\setminus
B(0,R)}|h(u_n)-h(u)-h(u_n-u)|^{\frac{2N}{N+\al}}|\phi|^{\frac{2N}{N+\al}}\le C\e\|\phi\|^{\frac{2N}{N+\al}}.
\ee
Finally, for $\e>0$ given above, there exists $C(\e)>0$ such that
\be\lab{bianhao2}
|h(t)|^{\frac{2N}{N+\al}}\le C(\e)|t|^{\frac{2N}{N+\al}}+\e|t|^{\frac{2N}{N+\al}\frac{2+\al}{N-2}},\,\, t\in\R.
\ee
Recalling that $u_n\rg u$ weakly in $H^1(\RN)$, up to a subsequence,
$u_n\rg u$ strongly in $L^{4N/(N+\al)}(B(0,R))$ and there exists $\omega\in
L^{4N/(N+\al)}(B(0,R))$ such that $|u_n(x)|,|u(x)|\le |\omega(x)|$ a.e. $x\in
B(0,R)$. Then it is easy to know for $n$ large,
\begin{align}\lab{bianhao4}
&\int_{B(0,R)}|h(u_n-u)|^{\frac{2N}{N+\al}}|\phi|^{\frac{2N}{N+\al}}\nonumber\\
&\le\int_{B(0,R)}\left(C(\e)|u_n-u|^{\frac{2N}{N+\al}}+\e|u_n-u|^{\frac{2N}{N+\al}\frac{2+\al}{N-2}}\right)|\phi|^{\frac{2N}{N+\al}}\le C\e\|\phi\|^{\frac{2N}{N+\al}}.
\end{align}
Moreover, let
$D_n:=\{x\in B(0,R): |u_n(x)-u(x)|\ge 1\}$,
then by \re{bianhao2},
{\allowdisplaybreaks
\begin{align*}
&\int_{D_n}|h(u_n)-h(u)|^{\frac{2N}{N+\al}}|\phi|^{\frac{2N}{N+\al}}\\
&\le\int_{D_n}\left[C(\e)(|u|^{\frac{2N}{N+\al}}+|u_n|^{\frac{2N}{N+\al}})+\e(|u_n|^{\frac{2N}{N+\al}\frac{2+\al}{N-2}}+|u|^{\frac{2N}{N+\al}\frac{2+\al}{N-2}})\right]|\phi|^{\frac{2N}{N+\al}}\\
&\le C\e\|\phi\|^{\frac{2N}{N+\al}}+2C(\e)\int_{D_n}|\omega|^{\frac{2N}{N+\al}}|\phi|^{\frac{2N}{N+\al}}\\
&\le
C\e\|\phi\|^{\frac{2N}{N+\al}}+2C(\e)\left(\int_{D_n}|\omega|^{\frac{4N}{N+\al}}\right)^{\frac{1}{2}}\left(\int_{\RN}|\phi|^{\frac{4N}{N+\al}}\right)^{\frac{1}{2}}.
\end{align*}}%
By $u_n\rg u$ a.e. $x\in B(0,R)$, we get that $|D_n|\rg 0$ as $n\rg
\iy$. Hence,
\be\lab{bianhao3}
\int_{D_n}|h(u_n)-h(u)|^{\frac{2N}{N+\al}}|\phi|^{\frac{2N}{N+\al}}\le C\e \|\phi\|^{\frac{2N}{N+\al}},\,\,\hbox{for $n$ large}.
\ee
On the other hand, for $\e$ given above, there exists $c(\e)>0$ such that
\begin{align*}
|h(u_n)-h(u)|^{\frac{2N}{N+\al}}&\le
\e(|u_n|^{\frac{2+\al}{N-2}\frac{2N}{N+\al}}+|u_n|^{\frac{2+\al}{N-2}\frac{2N}{N+\al}})\\
& \ \ \ +c(\e)(|u_n|^{\frac{2N}{N+\al}}+|u_n|^{\frac{2N}{N+\al}}).
\end{align*}
Noting that $|\{|u|\ge L\}|\rg 0$ as $L\rg \iy$, similar as above, there exists $L=L(\e)>0$ such that for all $n$,
$$
\int_{(B(0,R)\setminus D_n)\cap\{|u|\ge L\}}|h(u_n)-h(u)|^{\frac{2N}{N+\al}}|\phi|^{\frac{2N}{N+\al}}\le C\e\|\phi\|^{\frac{2N}{N+\al}}.
$$
By the Lebesgue dominated convergence theorem,
$$
\int_{(B(0,R)\setminus D_n)\cap\{|u|\le L\}}|h(u_n)-h(u)|^{\frac{2N}{N+\al}}|\phi|^{\frac{2N}{N+\al}}=o_n(1)\|\phi\|^{\frac{2N}{N+\al}},
$$
where $o_n(1)\rg0$ as $n\rg\iy$ uniformly for $\phi$. Then by \re{bianhao3},
$$
\int_{B(0,R)}|h(u_n)-h(u)|^{\frac{2N}{N+\al}}|\phi|^{\frac{2N}{N+\al}}\le C\e\|\phi\|^{\frac{2N}{N+\al}},\,\,\hbox{for $n$ large}.
$$
Then by \re{bianhao4} and for $n$ large,
$$
\int_{B(0,R)}|h(u_n)-h(u)-h(u_n-u)|^{\frac{2N}{N+\al}}|\phi|^{\frac{2N}{N+\al}}\le C\e\|\phi\|^{\frac{2N}{N+\al}},\,\,\hbox{for $n$ large}..
$$
Therefore, combing \re{bianhao5}, the proof is completed
\ep

\bl\lab{Asymptotic}
Let $\al\in(0,N)$ and $s\in(1,N/\al)$ and $\{g_n\}\in L^1(\RN)\cap L^s(\RN)$ be bounded both in $L^1(\RN)$ and $L^s(\RN)$ such that up to subsequences, for any bounded domain $\om\subset\RN$, $g_n\rg0$ strongly in $L^s(\om)$ as $n\rg\iy$, then passing to a subsequence if necessary, $(I_\al\ast g_n)(x)\rg0$ a. e. in $\RN$ as $n\rg\iy$.
\el
\bp It suffices to prove that for any fixed $k\in\mathbb{N}^+$, passing to a subsequence if necessary, $(I_\al\ast g_n)(x)\rg0$ a. e. in $B_k(0)$ as $n\rg\iy$. Let $k\in\mathbb{N}^+$ be fixed and for any $\dd>0$, there exists $K=K(\dd)>k$ such that
$$
A_\al\int_{\RN\setminus B_K(x)}\frac{|g_n(y)|}{|x-y|^{N-\al}}\,\ud y\le\dd,\,\,\mbox{for any}\,\, x\in\RN, n\in\mathbb{N}^+,
$$
Obviously, $B_K(x)\subset B_{2K}(0)$ for any $x\in B_K(0)$. Noting that $g_n\chi_{B_{2K}(0)}\in L^s(\RN)$, by Remark \ref{rr1},
$$
\|I_\al\ast(|g_n|\chi_{B_{2K}(0)})\|_{L^{\frac{Ns}{N-\al s}}(\RN)}\le C\|g_n\|_{L^s(B_{2K}(0))},
$$
where $C$ depends only on $N,\al$. It follows that up to a subsequence, $I_\al\ast(|g_n|\chi_{B_{2K}(0)})\rg0$ strongly in $L^{\frac{Ns}{N-\al s}}(\RN)$ and a. e. in $B_k(0)$ as $n\rg\iy$. Then for almost every $x\in B_k(0)$,
{\allowdisplaybreaks
\begin{align*}
&\limsup_{n\rg\iy}|(I_\al\ast g_n)(x)|\\
&\le A_\al\limsup_{n\rg\iy}\left(\int_{B_K(x)}\frac{|g_n(y)|}{|x-y|^{N-\al}}\,\ud y+\int_{\RN\setminus B_K(x)}\frac{|g_n(y)|}{|x-y|^{N-\al}}\,\ud y\right)\\
&\le\dd+A_\al\limsup_{n\rg\iy}\int_{B_K(x)}\frac{|g_n(y)|}{|x-y|^{N-\al}}\,\ud y\\
&\le\dd+A_\al\limsup_{n\rg\iy}\int_{B_{2K}(0)}\frac{|g_n(y)|}{|x-y|^{N-\al}}\,\ud y\\
&=\dd+\limsup_{n\rg\iy}[I_\al\ast(|g_n|\chi_{B_{2K}}(0))](x)=\dd.
\end{align*}}
Since $\dd$ is arbitrary, the proof is completed.
\ep
\noindent{\it Proof of Lemma \ref{Split}.}
Let
$$
f_1(t)=f(t)-|t|^{\frac{4+\al-N}{N-2}}t,\,\, F_1(t)=\int_0^tf_1(s)\,\ud s,\, t\in\R,
$$
then we observe that for any $\phi\in C_0^\iy(\RN)$,
\begin{align*}
\int_{\RN}[I_\al\ast F(u_n)]f(u_n)\phi=\int_{\RN}[I_\al\ast F(u_n)]f_1(u_n)\phi+\int_{\RN}[I_\al\ast F(u_n)]|u_n|^{\frac{4+\al-N}{N-2}}u_n\phi.
\end{align*}
{\bf Step 1.} We claim that
\begin{align*}
\int_{\RN}[I_\al\ast F(u_n)]|u_n|^{\frac{4+\al-N}{N-2}}u_n\phi&=\int_{\RN}[I_\al\ast F(u_n-u)]|u_n-u|^{\frac{4+\al-N}{N-2}}(u_n-u)\phi\\
&\,\,\,\,\,\,+\int_{\RN}[I_\al\ast F(u)]|u|^{\frac{4+\al-N}{N-2}}u\phi+o_n(1)\|\phi\|,
\end{align*}
where $o_n(1)\rg 0$ uniformly for any $\phi\in C_0^\iy(\RN)$ as $n\rg\iy$. Noting that $\al>N-4$, by $(ii)_1$ of Lemma \ref{Brezis} with $h(t)=f(t)$, $q=(2+\al)/(N-2), r=2N/(N-2)$,
\be\lab{yong}
\lim_{n\rg\iy}\int_{\RN}\big|F(u_n)-F(u_n-u)-F(u)\big|^{\frac{2N}{N+\al}}=0.
\ee
Then for $v_n=|u_n|^{\frac{4+\al-N}{N-2}}u_n, |u_n-u|^{\frac{4+\al-N}{N-2}}(u_n-u)$ or $|u|^{\frac{4+\al-N}{N-2}}u$, there exists $C>0$ such that
\begin{align*}
\int_{\RN}|v_n\phi|^{\frac{2N}{N+\al}}&\le\left(\int_{\RN}|v_n|^{\frac{2N}{2+\al}}\right)^{\frac{2+\al}{N+\al}}
\left(\int_{\RN}|\phi|^{\frac{2N}{N-2}}\right)^{\frac{N-2}{N+\al}}\\
&\le C\left(\int_{\RN}|\phi|^{\frac{2N}{N-2}}\right)^{\frac{N-2}{N+\al}}
\end{align*}
which follows that
{\allowdisplaybreaks
\begin{align*}
&\left|\int_{\RN}[I_\al\ast(F(u_n)-F(u_n-u)-F(u))]v_n\phi\right|\\
&\le C\left(\int_{\RN}\big|F(u_n)-F(u_n-u)-F(u))\big|^{\frac{2N}{N+\al}}\right)^{\frac{N+\al}{2N}}
\left(\int_{\RN}|v_n\phi|^{\frac{2N}{N+\al}}\right)^{\frac{N+\al}{2N}}\\
&=o_n(1)\left(\int_{\RN}|v_n\phi|^{\frac{2N}{N+\al}}\right)^{\frac{N+\al}{2N}}=o_n(1)\|\phi\|,
\end{align*}}%
where $o_n(1)\rg 0$ uniformly for any $\phi\in C_0^\iy(\RN)$ as $n\rg\iy$.
\vskip0.1in
On the other hand, by virtue of $(i)$ of Lemma \ref{Brezis} with $q=(2+\al)/(N-2)$ and $r=2N/(N-2)$,
$$
\lim_{n\rg\iy}\int_{\RN}\left||u_n|^{\frac{4+\al-N}{N-2}}u_n-|u_n-u|^{\frac{4+\al-N}{N-2}}(u_n-u)-|u|^{\frac{4+\al-N}{N-2}}u\right|^{\frac{2N}{2+\al}}=0.
$$
For $w_n=F(u_n), F(u_n-u)$ or $F(u)$, it is easy to know $\{w_n\}\subset L^{2N/(N+\al)}(\RN)$ is bounded in $L^{2N/(N+\al)}(\RN)$. By the Hardy-Littlewood-Sobolev inequality and H\"older's inequality, there exists $C>0$ such that
{\allowdisplaybreaks
\begin{align*}
&\left|\int_{\RN}[I_\al\ast w_n][|u_n|^{\frac{4+\al-N}{N-2}}u_n-|u_n-u|^{\frac{4+\al-N}{N-2}}(u_n-u)-|u|^{\frac{4+\al-N}{N-2}}u]\phi\right|\\
&\le C\left(\int_{\RN}\left||u_n|^{\frac{4+\al-N}{N-2}}u_n-|u_n-u|^{\frac{4+\al-N}{N-2}}(u_n-u)
-|u|^{\frac{4+\al-N}{N-2}}u\right|^{\frac{2N}{N+\al}}|\phi|^{\frac{2N}{N+\al}}\right)^{\frac{N+\al}{2N}}\\
&\le C\left(\int_{\RN}\left||u_n|^{\frac{4+\al-N}{N-2}}u_n-|u_n-u|^{\frac{4+\al-N}{N-2}}(u_n-u)
-|u|^{\frac{4+\al-N}{N-2}}u\right|^{\frac{2N}{2+\al}}\right)^{\frac{2+\al}{2N}}\left(\int_{\RN}|\phi|^{\frac{2N}{N-2}}\right)^{\frac{N-2}{2N}}\\
&=o_n(1)\|\phi\|,
\end{align*}}
where $o_n(1)\rg 0$ uniformly for any $\phi\in C_0^\iy(\RN)$ as $n\rg\iy$. Then we get
\begin{align*}
&\int_{\RN}[I_\al\ast F(u_n)]|u_n|^{\frac{4+\al-N}{N-2}}u_n\phi\\
&=\int_{\RN}[I_\al\ast F(u_n-u)]|u_n-u|^{\frac{4+\al-N}{N-2}}
(u_n-u)\phi+\int_{\RN}[I_\al\ast F(u)]|u|^{\frac{4+\al-N}{N-2}}u\phi\\
&\ \ \ +\int_{\RN}[I_\al\ast F(u_n-u)]|u|^{\frac{4+\al-N}{N-2}}u\phi
+\int_{\RN}[I_\al\ast F(u)]|u_n-u|^{\frac{4+\al-N}{N-2}}(u_n-u)\phi+o_n(1)\|\phi\|,
\end{align*}
where $o_n(1)\rg 0$ uniformly for any $\phi\in C_0^\iy(\RN)$ as $n\rg\iy$. Noting that $F(u)\in L^{2N/(N+\al)}(\RN)$, by Remark \ref{rr1},
$
|I_\al\ast F(u)|^{\frac{2N}{N+2}}\in L^{\frac{N+2}{N-\al}}(\RN).
$
By virtue of Lemma \ref{almost}, $|u_n-u|^{\frac{2N(2+\al)}{(N-2)(N+2)}}\rg 0$ weakly in $L^{(N+2)/(2+\al)}(\RN)$ as $n\rg0$. It follows that
$$
\lim_{n\rg\iy}\int_{\RN}|I_\al\ast F(u)|^{\frac{2N}{N+2}}|u_n-u|^{\frac{2N(2+\al)}{(N-2)(N+2)}}=0,
$$
which implies that
\begin{align*}
&\left|\int_{\RN}[I_\al\ast F(u)]|u_n-u|^{\frac{4+\al-N}{N-2}}(u_n-u)\phi\right|\\
&\le\left(\int_{\RN}|I_\al\ast F(u)|^{\frac{2N}{N+2}}|u_n-u|^{\frac{2N(2+\al)}{(N-2)(N+2)}}\right)^{\frac{N+2}{2N}}
\left(\int_{\RN}|\phi|^{\frac{2N}{N-2}}\right)^{\frac{N-2}{2N}}=o_n(1)\|\phi\|,
\end{align*}
where $o_n(1)\rg 0$ uniformly for any $\phi\in C_0^\iy(\RN)$ as $n\rg\iy$.
\vskip0.1in
Meanwhile, since $\al\in((N-4)_+,N)$, for $s\in(1, \frac{2N}{N+\al})\subset(1,\frac{N}{\al})$, by Rellich's theorem, up to sequences, for any bounded domain $\om\subset\RN$, $F(u_n-u)\rg0$ strongly in $L^s(\om)$ as $n\rg\iy$. By Lemma \ref{Asymptotic}, up to a sequence, $I_\al\ast F(u_n-u)\rg 0$ a. e. in $\RN$ as $n\rg0$. By Remark \ref{rr1},
$$
\sup_n\||I_\al\ast F(u_n-u)|^{\frac{2N}{N+2}}\|_{L^{\frac{N+2}{N-\al}}(\RN)}\le C\sup_n\|F(u_n-u)\|_{L^{\frac{2N}{N+\al}}(\RN)}<\iy,
$$
which yields by Lemma \ref{almost} that $|I_\al\ast F(u_n-u)|^{\frac{2N}{N+2}}\rg0$ weakly in $L^{\frac{N+2}{N-\al}}(\RN)$ as $n\rg\iy$. Noting that $|u|^{\frac{2+\al}{N-2}\frac{2N}{N+2}}\in L^{\frac{N+2}{2+\al}}(\RN)$,
\be\lab{guji1}
\lim_{n\rg\iy}\int_{\RN}|I_\al\ast F(u_n-u)|^{\frac{2N}{N+2}}|u|^{\frac{2+\al}{N-2}\frac{2N}{N+2}}=0,
\ee
and by H\"older's inequality,
\begin{align*}
&\left|\int_{\RN}[I_\al\ast F(u_n-u)]|u|^{\frac{4+\al-N}{N-2}}u\phi\right|\\
&\le\left(\int_{\RN}|I_\al\ast F(u_n-u)|^{\frac{2N}{N+2}}|u|^{\frac{2+\al}{N-2}\frac{2N}{N+2}}\right)^{\frac{N+2}{2N}}
\left(\int_{\RN}|\phi|^{\frac{2N}{N-2}}\right)^{\frac{N-2}{2N}}\\
&=o_n(1)\|\phi\|,
\end{align*}
where $o_n(1)\rg 0$ uniformly for any $\phi\in C_0^\iy(\RN)$ as $n\rg\iy$. The claim is concluded.
\vskip0.1in
{\bf Step 2.} We claim that
\begin{align}\lab{zh}
\int_{\RN}[I_\al\ast F(u_n)]f_1(u_n)\phi&=\int_{\RN}[I_\al\ast F(u_n-u)]f_1(u_n-u)\phi\\
&\ \ \ +\int_{\RN}[I_\al\ast F(u)]f_1(u)\phi+o_n(1)\|\phi\|,\nonumber
\end{align}
where $o_n(1)\rg 0$ uniformly for any $\phi\in C_0^\iy(\RN)$ as $n\rg\iy$. First, we prove that
\begin{align}\lab{yaoyong1}
\bcs
&\int_{\RN}[I_\al\ast(F(u_n)-F(u_n-u)-F(u))]f_1(u_n)\phi=o_n(1)\|\phi\|,\\
&\int_{\RN}[I_\al\ast(F(u_n)-F(u_n-u)-F(u))]f_1(u_n-u)\phi=o_n(1)\|\phi\|,\\
&\int_{\RN}[I_\al\ast(F(u_n)-F(u_n-u)-F(u))]f_1(u)\phi=o_n(1)\|\phi\|,
\ecs
\end{align}
where $o_n(1)\rg 0$ uniformly for any $\phi\in C_0^\iy(\RN)$ as $n\rg\iy$. We only prove the first quantity above. Other quantities can be proved in a similar way. Observe that there exists $\dd\in(0,1)$ and $C(\dd)>0$ such that $|f_1(t)|\le|t|$ for $|t|\le\dd$ and $|f_1(t)|\le C(\dd)|t|^{(2+\al)/(N-2)}$ for $|t|\ge\dd$. Noting that $\al\in((N-4)_+,N)$, we know $2<4N/(N+\al)<2N/(N-2)$. Then for any $\phi\in C_0^\iy(\RN)$, there exists $C>0$ (independent of $\phi, n$) such that
{\allowdisplaybreaks
\begin{align*}
&\int_{\RN}|f_1(u_n)\phi|^{\frac{2N}{N+\al}}=\int_{\{x\in\RN: |u_n(x)|\le\dd\}}|f_1(u_n)\phi|^{\frac{2N}{N+\al}}+\int_{\{x\in\RN: |u_n(x)|\ge\dd\}}|f_1(u_n)\phi|^{\frac{2N}{N+\al}}\\
&\le\int_{\{x\in\RN: |u_n(x)|\le\dd\}}|u_n\phi|^{\frac{2N}{N+\al}}+[C(\dd)]^{\frac{2N}{N+\al}}\int_{\{x\in\RN: |u_n(x)|\ge\dd\}}|u_n|^{\frac{2N(2+\al)}{(N-2)(N+\al)}}|\phi|^{\frac{2N}{N+\al}}\\
&\le\left(\int_{\RN}|u_n|^{\frac{4N}{N+\al}}\right)^{\frac{1}{2}}\left(\int_{\RN}|\phi|^{\frac{4N}{N+\al}}\right)^{\frac{1}{2}}
+[C(\dd)]^{\frac{2N}{N+\al}}\left(\int_{\RN}|u_n|^{\frac{2N}{N-2}}\right)^{\frac{2+\al}{N+\al}}\left(\int_{\RN}|\phi|^{\frac{2N}{N-2}}\right)^{\frac{N-2}{N+\al}}\\
&\le C\|\phi\|^{\frac{2N}{N+\al}},\,\,\mbox{for all}\,\, n=1,2,\cdots.
\end{align*}}%
It follows that
$$
\left(\int_{\RN}|f_1(u_n)\phi|^{\frac{2N}{N+\al}}\right)^{\frac{N+\al}{2N}}\le C\|\phi\|\,\,\,\mbox{uniformly for all}\,\,\phi\in C_0^\iy(\RN), n=1,2,\cdots.
$$
Then by the Hardy-Littlewood-Sobolev inequality and \re{yong},
\begin{align*}
&\left|\int_{\RN}[I_\al\ast(F(u_n)-F(u_n-u)-F(u))]f_1(u_n)\phi\right|\\
&\le C\left(\int_{\RN}\big|F(u_n)-F(u_n-u)-F(u)\big|^{\frac{2N}{N+\al}}\right)^{\frac{N+\al}{2N}}
\left(\int_{\RN}|f_1(u_n)\phi|^{\frac{2N}{N+\al}}\right)^{\frac{N+\al}{2N}}\\
&=o_n(1)\|\phi\|,
\end{align*}
where $o_n(1)\rg 0$ uniformly for any $\phi\in C_0^\iy(\RN)$ as $n\rg\iy$. So \re{yaoyong1} holds.

Second, we prove that
\begin{align}\lab{yaoyong2}
\bcs
&\int_{\RN}(I_\al\ast F(u_n))[f_1(u_n)-f_1(u_n-u)-f_1(u)]\phi=o_n(1)\|\phi\|,\\
&\int_{\RN}(I_\al\ast F(u_n-u))[f_1(u_n)-f_1(u_n-u)-f_1(u)]\phi=o_n(1)\|\phi\|,\\
&\int_{\RN}(I_\al\ast F(u))[f_1(u_n)-f_1(u_n-u)-f_1(u)]\phi=o_n(1)\|\phi\|,
\ecs
\end{align}
where $o_n(1)\rg 0$ uniformly for any $\phi\in C_0^\iy(\RN)$ as $n\rg\iy$. Similar as in Step 1, by the Hardy-Littlewood-Sobolev inequality and $(ii)_2$ of Lemma \ref{Brezis}, there exists $C>0$ such that
\begin{align*}
&\left|\int_{\RN}(I_\al\ast F(u_n))[f_1(u_n)-f_1(u_n-u)-f_1(u)]\phi\right|\\
&\le C\left(\int_{\RN}\left|f_1(u_n)-f_1(u_n-u)-f_1(u)\right|^{\frac{2N}{N+\al}}|\phi|^{\frac{2N}{N+\al}}\right)^{\frac{N+\al}{2N}}\\
&=o_n(1)\|\phi\|,
\end{align*}
where $o_n(1)\rg 0$ uniformly for any $\phi\in C_0^\iy(\RN)$ as $n\rg\iy$. So the first quantity of \re{yaoyong2} is concluded.

Then, combing \re{yaoyong1} and \re{yaoyong2}, we have
\begin{align*}
&\int_{\RN}[I_\al\ast F(u_n)]f_1(u_n)\phi\\
&=\int_{\RN}[I_\al\ast F(u_n-u)]f_1(u_n-u)\phi+\int_{\RN}[I_\al\ast F(u)]f_1(u)\phi\\
&\ \ \ +\int_{\RN}[I_\al\ast F(u_n-u)]f_1(u)\phi
+\int_{\RN}[I_\al\ast F(u)]f_1(u_n-u)\phi+o_n(1)\|\phi\|,
\end{align*}
where $o_n(1)\rg 0$ uniformly for any $\phi\in C_0^\iy(\RN)$ as $n\rg\iy$. To conclude the proof of \re{zh}, it suffices to prove
\be\lab{zh1}
\int_{\RN}[I_\al\ast F(u_n-u)]f_1(u)\phi=o_n(1)\|\phi\|,
\ee
and
\be\lab{zh2}
\int_{\RN}[I_\al\ast F(u)]f_1(u_n-u)\phi=o_n(1)\|\phi\|,
\ee
where $o_n(1)\rg 0$ uniformly for any $\phi\in C_0^\iy(\RN)$ as $n\rg\iy$. Notice that for any $\e\in(0,1)$, there exist $\dd\in(0,1)$ and $C_\e>0$ such that $|f_1(t)|\le\e|t|$ for $|t|\le\dd$ and $|f_1(t)|\le C_\e|t|^{(2+\al)/(N-2)}$ for $|t|\ge\dd$. Then for any $\phi\in C_0^\iy(\RN)$, by the Hardy-Littlewood-Sobolev inequality and H\"older's inequality,
{\allowdisplaybreaks
\begin{align*}
&\left|\int_{\RN}[I_\al\ast F(u_n-u)]f_1(u)\phi\right|\\
&\le\e\int_{\{x\in\RN:|u(x)|\le\dd\}}|I_\al\ast F(u_n-u)||u\phi|+C_\e\int_{\{x\in\RN:|u(x)|\ge\dd\}}|I_\al\ast F(u_n-u)||u|^{\frac{2+\al}{N-2}}|\phi|\\
&\le\e\|F(u_n-u)\|_{L^{\frac{2N}{N+\al}}(\RN)}\left(\int_{\{x\in\RN:|u(x)|\le\dd\}}|u\phi|^{\frac{2N}{N+\al}}\right)^{\frac{N+\al}{2N}}\\
&\ \ \ \ \ +C_\e\left(\int_{\RN}|I_\al\ast F(u_n-u)|^{\frac{2N}{N+2}}|u|^{\frac{2+\al}{N-2}\frac{2N}{N+2}}\right)^{\frac{N+2}{2N}}
\left(\int_{\RN}|\phi|^{\frac{2N}{N-2}}\right)^{\frac{N-2}{2N}}.
\end{align*}}%
Similar as above, there exists $c>0$ (independent of $\phi,\dd,\e$) such that
$$
\int_{\{x\in\RN:|u(x)|\le\dd\}}|u\phi|^{\frac{2N}{N+\al}}\le c\|\phi\|^{\frac{2N}{N+\al}}.
$$
Then by \re{guji1}, there exists $\ti{C}>0$ (independent of $\phi,\e$) such that
$$
\limsup_{n\rg\iy}\left|\int_{\RN}[I_\al\ast F(u_n-u)]f_1(u)\phi\right|\le\ti{C}\e\|\phi\|.
$$
It follows that
$$
\int_{\RN}[I_\al\ast F(u_n-u)]f_1(u)\phi=o_n(1)\|\phi\|,
$$
where $o_n(1)\rg 0$ uniformly for any $\phi\in C_0^\iy(\RN)$ as $n\rg\iy$. Similarly, \re{zh2} can be proved.
\vskip0.1in
Therefore, the proof of Lemma \ref{Split} is completed.
\qed

\subsection{Ground state solutions}
Since we are looking positive ground solutions of \re{lb1}, in this section, we may assume that $f$ is odd in $\RN$. We adapt the monotonicity trick due to L. Jeanjean \cite{Jeanjean} to seek ground state solutions of \re{lb1}.

For $\la\in[1/2,1]$, we consider the family of functionals as follows.
$$
I_\la(u)=\frac{1}{2}\int_{\RN}|\na u|^2+au^2-\frac{\la}{2}\int_{\RN} (I_\al\ast F(u))F(u),\ \ u\in H^1(\RN).
$$
Obviously, if $f$ satisfies the assumptions of Theorem \ref{Theorem 1}, for $\la\in[1/2,1]$ $I_\la\in C^1(H^1(\RN),\R)$ and every critical point of $I_\la$ is a weak solution of
\be\lab{q2a}
-\DD u+au=\la(I_\al\ast F(u))f(u),\ \ \ \ u\in H^1(\RN).
\ee
To guarantee the existence of critical points to $I_\la$, we recall the following abstract result, which was introduced by L. Jeanjean \cite{Jeanjean}.

\begin{theoremletter}\lab{Monotonicity}{\bf[see
\cite{Jeanjean}]} {\it Let $X$ be a Banach space equipped with a norm $\|\cdot\|_X$, $J\subset\R^+$ be an interval and a family of $C^1$-class functionals $\{I_\la\}_{\la\in J}$ on $X$ of the form
$$
I_\la(u)=A(u)-\la B(u),\,\, u\in X.
$$
Assume that $B(u)\ge0$ for any $u\in X$, one of $A, B$ is coercive in $X$ and there are two points $v_1,v_2\in X$ such that for any $\la\in J$,
$$
c_\la:=\inf_{\g\in\G}\max_{t\in[0,1]}I_\la(\g(t))>\max\{I_\la(v_1),I_\la(v_2\},
$$
where $\G:=\{\g\in C([0,1], X): \g(0)=v_1, \g(1)=v_2\}$. Then for almost every $\la\in J$, $I_\la$ admits a bounded Palais-Smale sequence for the level $c_\la$, Namely, there exists $\{v_n\}\subset X$ such that
\begin{itemize}
\item [{\rm(i)}] $\{v_n\}$ is bounded in $X$,
\item [{\rm(ii)}] $I_\la(v_n)\rg c_\la$ and $I_\la'(v_n)\rg0$ in $X^{-1}$ as $n\rg\iy$.
\end{itemize}
Moreover, $c_\la$ is continuous from the left-hand side with respect to $\la\in[1/2,1]$.
}
\end{theoremletter}
In the following, we use Theorem \ref{Monotonicity} to seek nontrivial weak solutions of \re{q2a} for almost every $\la\in [1/2,1]$. Then by passing to the limit, we get nontrivial solutions of the original problem \re{lb1}. In what follows, let $X=H^1(\RN)$ and
$$
A(u)=\frac{1}{2}\int_{\RN}|\na u|^2+au^2,\,\, B(u)=\frac{1}{2}\int_{\RN} (I_\al\ast F(u))F(u).
$$
Obviously, $A(u)\rg+\iy$ as $\|u\|\rg\iy$. Thanks to $(F3)$, $B(u)\ge0$ for any $u\in H^1(\RN)$. Moreover, by $(F1)$-$(F2)$, for any $\e>0$, there exists $C_\e>0$ such that $F(t)\le\e |t|^{(N+\al)/N}+C_\e|t|^{(N+\al)/(N-2)}$ for any $t\in\R$. Then similar as in \cite{MV1}, there exists $\dd>0$ such that
$$
\int_{\RN} (I_\al\ast F(u))F(u)\le\frac{1}{2}\|u\|^2,\,\, \mbox{if}\,\, \|u\|\le\dd,
$$
and therefore for any $u\in H^1(\RN)$ and $\la\in J$,
\be\lab{negative}
I_\la(u)\ge\frac{1}{4}\int_{\RN}|\na u|^2+au^2>0,\,\, \mbox{if}\,\, 0<\|u\|\le\dd.
\ee
On the other hand, taking a fixed $0\not=u_0\in H^1(\RN)$ and for any $\la\in J, t>0$, by $(F3)$,
$$
I_\la(\la u_0)\le\frac{t^2}{2}\int_{\RN}|\na u_0|^2+a|u_0|^2-\frac{t^{\frac{2N+2\al}{N-2}}}{4}\left(\frac{N-2}{N+\al}\right)^2\int_{\RN} (I_\al\ast |u_0|^{\frac{N+\al}{N-2}})|u_0|^{\frac{N+\al}{N-2}}
$$
and $I_\la(tu_0)\rg-\iy$ as $t\rg\iy$. Then there exists $t_0>0$ (independent of $\la$) such that $I_\la(t_0u_0)<0$, $\la\in J$ and $\|t_0u_0\|>\dd$. Let
$$
c_\la:=\inf_{\g\in\G}\max_{t\in[0,1]}I_\la(\g(t)),
$$
where $\G:=\{\g\in C([0,1], X): \g(0)=0, \g(1)=t_0u_0\}$.
\br
Here we remark that $c_\la$ is independent of $u_0$. In fact, let
$$
d_\la:=\inf_{\g\in\G_1}\max_{t\in[0,1]}I_\la(\g(t)),
$$
where $\G_1:=\{\g\in C([0,1], X): \g(0)=0, I_\la(\g(1))<0\}$. Obviously, $d_\la\le c_\la$. On the other hand, for any $\g\in\G_1$, it follows from \re{negative} that $\|\g(1)\|>\dd$. Due to the path connectedness of $H^1(\RN)$, there exists $\ti{\g}\in C([0,1], H^1(\RN))$ such that $\ti{\g}(t)=\g(2t)$ if $t\in[0,1/2]$, $\|\ti{\g}(t)\|>\dd$ if $t\in[1/2,1]$ and $\ti{\g}(1)=t_0u_0$. Then $\ti{\g}\in\G$ and $\max_{t\in[0,1]}I_\la(\ti{\g}(t))=\max_{t\in[0,1]}I_\la(\g(t))$, which implies that $c_\la\le d_\la$ and so $d_\la=c_\la$ for any $\la\in J$.
\er

By \re{negative}, $c_\la>\dd^2/4$ for any $\la\in J$. Then, as a consequence of Theorem \ref{Monotonicity}, we have\
\bl\lab{bps}
Assume $(F1)$-$(F3)$, for almost every $\la\in J=[1/2,1]$, problem \re{q2a} possesses a bounded Palais-Smale sequence at the level $c_\la$. Namely, there exists $\{u_n\}\subset H^1(\RN)$ such that
\begin{itemize}
\item [{\rm(i)}]  $\{u_n\}$ is bounded in $H^1(\RN)$,
\item [{\rm(ii)}] $I_\la(u_n)\rg c_\la$ and $I_\la'(u_n)\rg0$ in $H^{-1}(\RN)$ as $n\rg\iy$.
\end{itemize}
\el
In the following, in the spirit of \cite{JT, Lions1}, we give a decomposition of $\{u_n\}$ above, which plays a crucial role in the existence of ground states to \re{lb1}. But due to the present of a nonlocal and critical nonlinearity in the Hardy-Littlewood-Sobolev case, in contrast with the local term in the subcritical case, the proof becomes much more complicated.
\bo\lab{describe}
With the same assumptions in Theorem \ref{Theorem 2} and let $\la\in[1/2,1]$, $\{u_n\}$ be given in Lemma \ref{bps}. Assume that $u_n\rg u_\la$ weakly but not strongly in $H^1(\RN)$ as $n\rg\iy$, then up to a sequence, there exist $k\in\mathbb{N}^+$, $\{x_n^j\}_{j=1}^k\subset\RN$ and $\{v_\la^j\}_{j=1}^k\subset H^1(\RN)$ such that
\begin{itemize}
\item [{\rm(i)}]  $I_\la'(u_\la)=0$ in $H^{-1}(\RN)$,
\item [{\rm(ii)}] $v_\la^j\not=0$ and $I_\la'(v_\la^j)=0$ in $H^{-1}(\RN)$, $j=1,2,\cdots,k$,
\item [{\rm(iii)}] $c_\la=I_\la(u_\la)+\sum_{j=1}^kI_\la(v_\la^j)$,
\item [{\rm(iv)}] $\left\|u_n-u_\la-\sum_{j=1}^kv_\la^j(\cdot-x_n^j)\right\|\rg0$ as $n\rg\iy$,
\item [{\rm(v)}] $|x_n^j|\rg\iy$ and $|x_n^i-x_n^j|\rg\iy$ as $n\rg\iy$ for any $i\not=j$.
\end{itemize}
\eo
\noindent Before proving Proposition \ref{describe}, we give some lemmas.

\bl\lab{Pohozaev}
Let $\la\in[1/2,1]$ and $u_\la$ be any nontrivial weak solution of \re{q2a}, then $u_\la$ satisfies the following Pohoz\v{a}ev identity
$$
 \frac{N-2}{2}\int_{\mathbb{R}^N}|\na u_\la|^2+\frac{N}{2}a\int_{\mathbb{R}^N}|u_\la|^2=\frac{N+\al}{2}\la\int_{\mathbb{R}^N}(I_\al\ast F(u_\la))F(u_\la).
$$
Moreover, there exists $\beta,\g>0$ (independent of $\la\in[1/2,1]$) such that $\|u_\la\|\ge\beta$ and $I_\la(u_\la)\ge\g$ for any nontrivial solution $u_\la$, $\la\in[1/2,1]$.
\el
\bp
The proof of Pohoz\v{a}ev's identity is similar as in \cite[Theorem 3]{MV1}. We omit the detail here. Let $\la\in[1/2,1]$ and $u_\la$ be any nontrivial weak solution of \re{q2a}, then
\be\lab{bjiao}
\int_{\RN}|\na u_\la|^2+a|u_\la|^2\le\int_{\mathbb{R}^N}(I_\al\ast F(u_\la))f(u_\la)u_\la.
\ee
Thanks to $(F1)$-$(F2)$, for any $\e>0$, there exists $C_\e>0$ such that $F(t), tf(t)\le\e |t|^{(N+\al)/N}+C_\e|t|^{(N+\al)/(N-2)}$ for any $t\in\R$. Then similar as in \cite{MV1}, there exists $\beta>0$ such that
$$
\int_{\RN} (I_\al\ast F(u))f(u)u\le\frac{\|u\|^2}{2},\,\, \mbox{if}\,\, \|u\|\le\beta,
$$
which yields by \re{bjiao} that $\|u_\la\|\ge\beta$. Meanwhile, by Pohoz\v{a}ev's identity,
$$
I_\la(u_\la)=\frac{2+\al}{2(N+\al)}\int_{\RN}|\na u_\la|^2+\frac{\al a}{2(N+\al)}\int_{\RN}|u_\la|^2
$$
and then we conclude the proof.
\ep
Let $\al\in(0,N)$. For any $u\in D^{1,2}(\RN)$, by the Hardy-Littlewood-Sobolev inequality and Sobolev's inequality,
\begin{align*}
\int_{\RN}(I_\al\ast|u|^{\frac{N+\al}{N-2}})|u|^{\frac{N+\al}{N-2}}&\le A_\al\mathcal{C}_\al\left(\int_{\RN}|u|^{\frac{2N}{N-2}}\right)^{\frac{N+\al}{N}}\\
&\le A_\al\mathcal{C}_\al\mathcal{S}^{-\frac{N+\al}{N-2}}\left(\int_{\RN}|\na u|^2\right)^{\frac{N+\al}{N-2}},
\end{align*}
where
$$
\mathcal{S}:=\inf_{0\not=u\in D^{1,2}(\RN)}\frac{\int_{\RN}|\na u|^2}
{\left(\int_{\RN}|u|^{\frac{2N}{N-2}}\right)^{\frac{N-2}{N}}}.
$$
Then
$$
\mathcal{S}_\al:=\inf_{0\not=u\in D^{1,2}(\RN)}\frac{\int_{\RN}|\na u|^2}{\left(\int_{\RN}(I_\al\ast|u|^{\frac{N+\al}{N-2}})|u|^{\frac{N+\al}{N-2}}\right)^{\frac{N-2}{N+\al}}}
\ge\frac{\mathcal{S}}{(A_\al\mathcal{C}_\al)^{\frac{N-2}{N+\al}}}.
$$
In \cite[Lemma 1.2]{GY}, F. Gao and M. Yang proved that $\mathcal{S}_\al=\mathcal{S}/(A_\al\mathcal{C}_\al)^{\frac{N-2}{N+\al}}$ and can be achieved by the instanton
$$
U(x)=\frac{[N(N-2)]^{\frac{N-2}{4}}}{(1+|x|^2)^{\frac{N-2}{2}}}.
$$
Now, we give an upper estimate of $c_\la$.
\bl\lab{upper}
For any $\la\in[1/2,1]$, $\al\in(0,N)$ and assume
$$
q>\max\left\{1+\frac{\al}{N-2},\frac{N+\al}{2(N-2)}\right\},
$$
then we have
$$
c_\la<\frac{2+\al}{2(N+\al)}\left(\frac{N+\al}{N-2}\right)^{\frac{N-2}{2+\al}}\la^{\frac{2-N}{2+\al}}\mathcal{S}_\al^{\frac{N+\al}{2+\al}}.
$$
\el
\bp
Let $\varphi\in C_0^{\infty}(\RN)$ is a cut-off function with support $B_2$ such that $\varphi \equiv 1$ on $B_1$ and $0\le \varphi \le 1$ on $B_2$, where $B_r$ denotes the ball in $\RN$ of center at origin and radius $r$. Given $\e>0$, we set $\psi_\e(x)=\varphi(x)U_{\e}(x)$, where
$$
U_{\e}(x)=\frac{\big(N(N-2)\e^2\big)^{\frac{N-2}{4}}}{\big(\e^2+|x|^2\big)^{\frac{N-2}{2}}}.
$$
By \cite{Brezis2}(see also\cite[Lemma 1.46]{Willem1}), we have the following estimates:
\begin{align*}
\int_{\RN}|\nabla \psi_{\e}|^2=\mathcal{S}^{\frac{N}{2}}+
\bcs
O(\e^{N-2}), & \hbox{if $N\ge4$,}\\
K_1\e+O(\e^3), & \hbox{if $N=3$,}\\
\ecs
\end{align*}
$$
\int_{\RN}|\psi_{\e}|^{\frac{2N}{N-2}}=\mathcal{S}^{\frac{N}{2}}+O(\e^N),\,\,\, \hbox{if $N\ge3$,}
$$
and
\begin{align*}
\int_{\RN}|\psi_\e|^2=
\bcs
K_2\e^2+O(\e^{N-2}), & \hbox{if $N\ge5$,}\\
K_2\e^2|\ln\e|+O(\e^2), & \hbox{if $N=4$,}\\
K_2\e+O(\e^2), & \hbox{if $N=3$,}
\ecs
\end{align*}
where $K_1, K_2>0$. Then we get
\begin{align}\lab{esti}
\int_{\RN}|\na\psi_\e|^2+a|\psi_\e|^2=\mathcal{S}^{\frac{N}{2}}+
\bcs
aK_2\e^2+O(\e^{N-2}), & \hbox{if $N\ge5$,}\\
aK_2\e^2|\ln\e|+O(\e^2), & \hbox{if $N=4$,}\\
(K_1+aK_2)\e+O(\e^2), & \hbox{if $N=3$.}
\ecs
\end{align}
By direct computation, we know
$$
\left(\int_{\RN}|\psi_\e|^{\frac{2Nq}{N+\al}}\right)^{\frac{N+\al}{N}}=K_3\varepsilon^{N+\al-(N-2)q}+o(\varepsilon^{N+\al-(N-2)q}),
$$
and then by the Hardy-Littlewood-Sobolev inequality,
\begin{align}\lab{estimate5}
\int_{\RN}(I_\al\ast|\psi_\e|^{\frac{N+\al}{N-2}})|\psi_\e|^q&\le\mathcal{C}_\al\left(\int_{\RN}|\psi_{\e}|^{\frac{2N}{N-2}}\right)^{\frac{N+\al}{2N}}
\left(\int_{\RN}|\psi_{\e}|^{\frac{2Nq}{N+\al}}\right)^{\frac{N+\al}{2N}}\nonumber\\
&\le K_4\varepsilon^{\frac{N+\al-(N-2)q}{2}}+o(\varepsilon^{\frac{N+\al-(N-2)q}{2}}),
\end{align}
where $K_3,K_4>0$. Moreover, similar as in \cite{GY, GY1}, by direct computation, for some $K_5>0$,
\be\lab{estimate4}
\int_{\RN}(I_\al\ast|\psi_\e|^{\frac{N+\al}{N-2}})|\psi_\e|^{\frac{N+\al}{N-2}}
\ge(A_\al\mathcal{C}_\al)^{\frac{N}{2}}\mathcal{S}_\al^{\frac{N+\al}{2}}-K_5\varepsilon^{\frac{N+\al}{2}}+o(\varepsilon^{\frac{N+\al}{2}}).
\ee
Meanwhile,
\begin{align*}
&\int_{\RN}(I_\al\ast|\psi_\e|^{\frac{N+\al}{N-2}})|\psi_\e|^q\\
&\ge A_\al\left(\int_{\RN}\int_{\RN}\frac{U_\e^{\frac{N+\al}{N-2}}(x)U_\e^q(y)}{|x-y|^{N-\al}}\,\ud x\,\ud y-\int_{\RN\setminus B_1}\int_{B_1}\frac{U_\e^{\frac{N+\al}{N-2}}(x)U_\e^q(y)}{|x-y|^{N-\al}}\,\ud x\,\ud y\right.\\
&\,\,\,\,\,\,\,\,\,\,\,\,\,\,\,\,\,\,\,\,\,\,\,-\left.\int_{B_1}\int_{\RN\setminus B_1}\frac{U_\e^{\frac{N+\al}{N-2}}(x)U_\e^q(y)}{|x-y|^{N-\al}}\,\ud x\,\ud y-\int_{\RN\setminus B_1}\int_{\RN\setminus B_1}\frac{U_\e^{\frac{N+\al}{N-2}}(x)U_\e^q(y)}{|x-y|^{N-\al}}\,\ud x\,\ud y\right),
\end{align*}
where for some $\ti{K}_i>0$, $i=1,2,3,4$,
\begin{align*}
\bcs
\int_{\RN}\int_{\RN}\frac{U_\e^{\frac{N+\al}{N-2}}(x)U_\e^q(y)}{|x-y|^{N-\al}}\,\ud x\,\ud y=\ti{K}_1\e^{\frac{N+\al-(N-2)q}{2}},\\
\int_{\RN\setminus B_1}\int_{B_1}\frac{U_\e^{\frac{N+\al}{N-2}}(x)U_\e^q(y)}{|x-y|^{N-\al}}\,\ud x\,\ud y\le\ti{K}_2\e^{N+\al-\frac{N-2}{2}q}+o(\e^{N+\al-\frac{N-2}{2}q}),\\
\int_{B_1}\int_{\RN\setminus B_1}\frac{U_\e^{\frac{N+\al}{N-2}}(x)U_\e^q(y)}{|x-y|^{N-\al}}\,\ud x\,\ud y\le\ti{K}_3\e^{\frac{N-2}{2}q}+o(\e^{\frac{N-2}{2}q}),\\
\int_{\RN\setminus B_1}\int_{\RN\setminus B_1}\frac{U_\e^{\frac{N+\al}{N-2}}(x)U_\e^q(y)}{|x-y|^{N-\al}}\,\ud x\,\ud y\le\ti{K}_4\e^{\frac{N+\al+(N-2)q}{2}}+o(\e^{\frac{N+\al+(N-2)q}{2}}),
\ecs
\end{align*}
then for some $K_6>0$, we have
\be\lab{estimate41}
\int_{\RN}(I_\al\ast|\psi_\e|^{\frac{N+\al}{N-2}})|\psi_\e|^q\ge K_6\varepsilon^{\frac{N+\al-(N-2)q}{2}}+o(\varepsilon^{\frac{N+\al-(N-2)q}{2}}).
\ee
Here, we used the fact that $q>(N+\al)/[2(N-2)]$. Then for any $t>0$,
\begin{align*}
I_\la(t\psi_\e)\le&\frac{t^2}{2}\int_{\RN}|\na\psi_\e|^2+a|\psi_\e|^2-\frac{\mu\la}{q}\frac{N-2}{N+\al}t^{q+\frac{N+\al}{N-2}}\int_{\RN}(I_\al\ast\psi_\e^{\frac{N+\al}{N-2}})\psi_\e^q\\
&-\frac{t^{\frac{2(N+\al)}{N-2}}}{2}\left(\frac{N-2}{N+\al}\right)^2\la\int_{\RN}(I_\al\ast\psi_\e^{\frac{N+\al}{N-2}})\psi_\e^{\frac{N+\al}{N-2}}\\
:=&g_\e(t).
\end{align*}
Obviously, $g_\e(t)\rg-\iy$ as $t\rg+\iy$ and $g_\e(t)>0$ for $t>0$ small. Similar to \cite[Lemma 3.3]{Ruiz}, $g_\e$ has a unique critical point $t_\e$ in $(0,+\iy)$, which is the maximum point of $g_\e$. Meanwhile, by $g_\e'(t_\e)=0$,
\begin{align}\lab{estimate6}
&t_\e\int_{\RN}|\na\psi_\e|^2+a|\psi_\e|^2-\left(q+\frac{N+\al}{N-2}\right)\frac{\mu\la}{q}\frac{N-2}{N+\al}t_\e^{q+\frac{N+\al}{N-2}-1}\int_{\RN}(I_\al\ast\psi_\e^{\frac{N+\al}{N-2}})\psi_\e^q\nonumber\\
=&t_\e^{\frac{2(N+\al)}{N-2}-1}\frac{N-2}{N+\al}\la\int_{\RN}(I_\al\ast\psi_\e^{\frac{N+\al}{N-2}})\psi_\e^{\frac{N+\al}{N-2}}.
\end{align}
{\bf Claim.} There exists $t_0, t_1>0$(independent of $\e$) such that $t_\e\in[t_0,t_1]$ for $\e>0$ small. First, if $t_\e\rg0$ as $\e\rg0$, then by \re{esti}, \re{estimate5} and \re{estimate4}, there exist $c_1,c_2>0$(independent of $\e$) such that for $\e$ small,
$$
c_1t_\e\le c_2\varepsilon^{\frac{N+\al-(N-2)q}{2}}t_\e^{q+\frac{N+\al}{N-2}-1}+t_\e^{q+\frac{N+\al}{N-2}-1}\le 2t_\e^{q+\frac{N+\al}{N-2}-1},
$$
where we used the fact that $q<(N+\al)/(N-2)$. Then we get a contradiction. Second, by \re{estimate6},
\begin{align*}
\int_{\RN}|\na\psi_\e|^2+a|\psi_\e|^2\ge t_\e^{\frac{2(N+\al)}{N-2}-2}\frac{N-2}{N+\al}\la\int_{\RN}(I_\al\ast\psi_\e^{\frac{N+\al}{N-2}})\psi_\e^{\frac{N+\al}{N-2}},
\end{align*}
which implies combing \re{esti} and \re{estimate4} that $t_\e\le t_1$ for some $t_1>0$ and $\e$ small.

Then, by {\bf Claim} and \re{estimate41}, for some $K_7>0$,
$$
\frac{\mu\la}{q}\frac{N-2}{N+\al}t_\e^{q+\frac{N+\al}{N-2}}\int_{\RN}(I_\al\ast\psi_\e^{\frac{N+\al}{N-2}})\psi_\e^q\ge K_7\varepsilon^{\frac{N+\al-(N-2)q}{2}}+o(\varepsilon^{\frac{N+\al-(N-2)q}{2}}),
$$
which follows that
{\allowdisplaybreaks
\begin{align*}
&\max_{t\ge0}I_\la(t\psi_\e)=g_\e(t_\e)\\
&\le\frac{t_\e^2}{2}\int_{\RN}|\na\psi_\e|^2+a|\psi_\e|^2-K_7\varepsilon^{\frac{N+\al-(N-2)q}{2}}\\
&\,\,\,\,\,\,-\frac{t_\e^{\frac{2(N+\al)}{N-2}}}{2}\left(\frac{N-2}{N+\al}\right)^2\la\int_{\RN}(I_\al\ast\psi_\e^{\frac{N+\al}{N-2}})\psi_\e^{\frac{N+\al}{N-2}}
+o(\varepsilon^{\frac{N+\al-(N-2)q}{2}})\\
&\le\max_{t\ge0}\left[\frac{t^2}{2}\int_{\RN}|\na\psi_\e|^2+a|\psi_\e|^2-\frac{t^{\frac{2(N+\al)}{N-2}}}{2}\left(\frac{N-2}{N+\al}\right)^2\la\int_{\RN}(I_\al\ast\psi_\e^{\frac{N+\al}{N-2}})\psi_\e^{\frac{N+\al}{N-2}}\right]\\
&\,\,\,\,\,\,-K_7\varepsilon^{\frac{N+\al-(N-2)q}{2}}+o(\varepsilon^{\frac{N+\al-(N-2)q}{2}})\\
&=\frac{2+\al}{2(N+\al)}\left(\frac{N+\al}{N-2}\right)^{\frac{N-2}{2+\al}}\la^{^{\frac{2-N}{2+\al}}}\frac{\left(\int_{\RN}|\na\psi_\e|^2+a|\psi_\e|^2\right)^{\frac{N+\al}{2+\al}}}
{\left(\int_{\RN}(I_\al\ast\psi_\e^{\frac{N+\al}{N-2}})\psi_\e^{\frac{N+\al}{N-2}}\right)^{\frac{N-2}{2+\al}}}\\
&\,\,\,\,\,\,-K_7\varepsilon^{\frac{N+\al-(N-2)q}{2}}+o(\varepsilon^{\frac{N+\al-(N-2)q}{2}}).
\end{align*}
}%
On the other hand, by \re{esti} and \re{estimate4}, for some $K_8>0$,
{\allowdisplaybreaks
\begin{align*}
&\frac{\left(\int_{\RN}|\na\psi_\e|^2+a|\psi_\e|^2\right)^{\frac{N+\al}{2+\al}}}
{\left(\int_{\RN}(I_\al\ast\psi_\e^{\frac{N+\al}{N-2}})\psi_\e^{\frac{N+\al}{N-2}}\right)^{\frac{N-2}{2+\al}}}\\
&\le\mathcal{S}_\al^{\frac{N+\al}{2+\al}}+\bcs
K_8\e^{\min\{2,\frac{N+\al}{2}\}}+o(\varepsilon^{\min\{2,\frac{N+\al}{2}\}}),\,&\hbox{if $N\ge5$},\\
K_8\e^2|\ln\e|+o(\e^2|\ln\e|),\,&\hbox{if $N=4$},\\
K_8\e+o(\e),\,&\hbox{if $N=3$}.
\ecs
\end{align*}}%
Then, for some $K_9,K_{10}>0$,
{\allowdisplaybreaks
\begin{align*}
\max_{t\ge0}I_\la(t\psi_\e)&\le \frac{2+\al}{2(N+\al)}\left(\frac{N+\al}{N-2}\right)^{\frac{N-2}{2+\al}}\la^{^{\frac{2-N}{2+\al}}}\mathcal{S}_\al^{\frac{N+\al}{2+\al}}\\
&\,\,\,\,\,+\bcs
K_9\e^{\min\{2,\frac{N+\al}{2}\}}-K_{10}\varepsilon^{\frac{N+\al-(N-2)q}{2}}+o(\varepsilon^{\frac{N+\al-(N-2)q}{2}}),\,&\hbox{if $N\ge5$},\\
K_9\e^2|\ln\e|-K_{10}\varepsilon^{\frac{N+\al-(N-2)q}{2}}+o(\varepsilon^{\frac{N+\al-(N-2)q}{2}}),\,&\hbox{if $N=4$},\\
K_9\e-K_{10}\varepsilon^{\frac{N+\al-(N-2)q}{2}}+o(\varepsilon^{\frac{N+\al-(N-2)q}{2}}),\,&\hbox{if $N=3$},
\ecs\\
&<\frac{2+\al}{2(N+\al)}\left(\frac{N+\al}{N-2}\right)^{\frac{N-2}{2+\al}}\la^{^{\frac{2-N}{2+\al}}}\mathcal{S}_\al^{\frac{N+\al}{2+\al}},\,\,\hbox{if $\e>0$ small enough},
\end{align*}}%
where we used the fact that $N+\al-(N-2)q<\min\{2,\frac{N+\al}{2}\}$. Therefore, for any $\la\in[1/2,1]$ and $\e>0$ small enough, we get
$$
c_\la\le\max_{t\ge0}I_\la(t\psi_\e)<\frac{2+\al}{2(N+\al)}\left(\frac{N+\al}{N-2}\right)^{\frac{N-2}{2+\al}}\la^{^{\frac{2-N}{2+\al}}}\mathcal{S}_\al^{\frac{N+\al}{2+\al}}.
$$
The proof is completed.
\ep

\noindent {\it Proof of Proposition \ref{describe}.} Let $\la\in[1/2,1]$ and assume $u_n\rg u_\la$ weakly in $H^1(\RN)$ as $n\rg\iy$ but not strongly in $H^1(\RN)$ and satisfies $I_\la(u_n)\rg c_\la$ and $I_\la'(u_n)\rg0$ in $H^{-1}(\RN)$ as $n\rg\iy$.

{\bf Step 1.} We claim that $I_\la'(u_\la)=0$ in $H^{-1}(\RN)$. As a consequence of Lemma \ref{Split}, it suffices to show for any fixed $\phi\in C_0^\iy(\RN)$, up to a sequence,
$$
\int_{\RN}[I_\al\ast F(u_n-u)]f(u_n-u)\phi\rg0,\hbox{as $n\rg\iy$}.
$$
In fact, by $(F1)$-$(F2)$, there exists $C>0$ such that
$$
|f(t)|^{\frac{2N}{N+\al}}\le C(|t|^{\frac{2N}{N+\al}}+|t|^{\frac{2+\al}{N-2}\frac{2N}{N+\al}}),\,\ t\in\R.
$$
By virtue of the Hardy-Littlewood-Sobolev inequality and Rellich's theorem, up to a sequence, for some $C$(independent of $n$) such that
$$
\left|\int_{\RN}[I_\al\ast F(u_n-u)]f(u_n-u)\phi\right|\le C\left(\int_{\RN}|f(u_n-u)\phi|^{\frac{2N}{N+\al}}\right)^{\frac{N+\al}{2N}}\rg0,\hbox{as $n\rg\iy$}.
$$

{\bf Step 2.} Set $v_n^1:=u_n-u_\la$, we claim that
\be\lab{nonvanish}
\lim_{n\rg\iy}\sup_{z\in\RN}\int_{B_1(z)}|v_n^1|^2>0.
\ee
Otherwise, by Lions' lemma \cite[Lemma I.1]{Lions1}, $v_n^1\rg0$ strongly in $L^t(\RN)$ as $n\rg\iy$ for any $t\in(2,2N/(N-2))$. Noting that
$\lan I_\la'(u_n), v_n^1\ran\rg0$ as $n\rg\iy$ and $\lan I_\la'(u_\la), v_n^1\ran=0$ for any $n$, by virtue of Lemma \ref{Lieb} and Lemma \ref{Split}, we get
\be\lab{decom}
\bcs
c_\la=I_\la(u_\la)+I_\la(v_n^1)+o_n(1),\\
\|v_n^1\|^2=\la\int_{\RN}[I_\al\ast F(v_n^1)]f(v_n^1)v_n^1+o_n(1),
\ecs
\ee
where $o_n(1)\rg0$ as $n\rg\iy$. Now, we show that
$$
\lim_{n\rg\iy}\int_{\RN}[I_\al\ast F_1(v_n^1)]F_1(v_n^1)=0,
$$
where
$$
f_1(t)=f(t)-|t|^{\frac{4+\al-N}{N-2}}t,\,\, F_1(t)=\int_0^tf_1(s)\,\ud s,\, t\in\R.
$$
Notice that $4N/(N+\al)\in(2,2N/(N-2))$ and $f_1(t)=o(t)$ as $|t|\rg0$, $\lim_{|t|\rg\iy}|f_1(t)|/|t|^{\frac{2+\al}{N-2}}=0$. It is easy to know
$$
\lim_{n\rg\iy}\int_{\RN}|F_1(v_n^1)|^{\frac{2N}{N+\al}}=0,
$$
which yields by the Hardy-Littlewood-Sobolev inequality that there exists some $C>0$(independent of $n$), such that
$$
\left|\int_{\RN}[I_\al\ast F_1(v_n^1)]F_1(v_n^1)\right|\le C\left(\int_{\RN}|F_1(v_n^1)|^{\frac{2N}{N+\al}}\right)^{\frac{N+\al}{2N}}\rg0,\,\,\hbox{as $n\rg\iy$}.
$$
Similarly,
\begin{align*}
\bcs
\lim_{n\rg\iy}\int_{\RN}[I_\al\ast F_1(v_n^1)]|v_n^1|^{\frac{N+\al}{N-2}}=0,\\
\lim_{n\rg\iy}\int_{\RN}[I_\al\ast F_1(v_n^1)]f_1(v_n^1)v_n^1=0.
\ecs
\end{align*}
Then by \re{decom}, we get
\be\lab{decom1}
\bcs
c_\la=I_\la(u_\la)+\frac{1}{2}\|v_n^1\|^2-\frac{\la}{2}\left(\frac{N-2}{N+\al}\right)^2\int_{\RN}[I_\al\ast |v_n^1|^{\frac{N+\al}{N-2}}]|v_n^1|^{\frac{N+\al}{N-2}}+o_n(1),\\
\|v_n^1\|^2=\la\frac{N-2}{N+\al}\int_{\RN}[I_\al\ast |v_n^1|^{\frac{N+\al}{N-2}}]|v_n^1|^{\frac{N+\al}{N-2}}+o_n(1),
\ecs
\ee
where $o_n(1)\rg0$ as $n\rg\iy$. Recalling that $v_n^1\not\rg0$ strongly in $H^1(\RN)$ as $n\rg\iy$, let
$$
\lim_{n\rg\iy}\|v_n^1\|^2=\la\frac{N-2}{N+\al}\lim_{n\rg\iy}\int_{\RN}[I_\al\ast |v_n^1|^{\frac{N+\al}{N-2}}]|v_n^1|^{\frac{N+\al}{N-2}}=b,
$$
then $b>0$. Noting that
$$
\int_{\RN}|\na v_n^1|^2\ge\mathcal{S}_\al\left(\int_{\RN}[I_\al\ast|v_n^1|^{\frac{N+\al}{N-2}}]|v_n^1|^{\frac{N+\al}{N-2}}\right)^{\frac{N-2}{N+\al}},
$$
we know
$$
b\ge\left(\frac{N+\al}{N-2}\right)^{\frac{N-2}{2+\al}}\la^{\frac{2-N}{2+\al}}\mathcal{S}_\al^{\frac{N+\al}{2+\al}}.
$$
By Lemma \ref{Pohozaev} and \re{decom1}
$$
c_\la\ge\frac{2+\al}{2(N+\al)}\left(\frac{N+\al}{N-2}\right)^{\frac{N-2}{2+\al}}\la^{\frac{2-N}{2+\al}}\mathcal{S}_\al^{\frac{N+\al}{2+\al}},
$$
which is a contradiction. Thus, \re{nonvanish} is true.
\vskip0.1in
{\bf Step 3.} By \re{nonvanish} and $v_n^1\rg0$ weakly in $H^1(\RN)$ as $n\rg\iy$, there exists $\{z_n^1\}\subset\RN$ such that $|z_n^1|\rg\iy$ as $n\rg\iy$ and
$
\lim_{n\rg\iy}\int_{B_1(z_n^1)}|v_n^1|^2>0.
$
Let $u_n^1=v_n^1(\cdot+z_n^1)$, then up to a sequence, $u_n^1\rg v_\la^1$ weakly in $H^1(\RN)$ as $n\rg\iy$ for some $v_\la^1\not=0$. By virtue of Lemma \ref{Lieb} and Lemma \ref{Split}, it is easy to know
$$
I_\la(u_n^1)\rg c_\la-I_\la(u_\la),\,\, I_\la'(u_n^1)\rg0\,\,\hbox{in $H^{-1}(\RN)$ as $n\rg\iy$}.
$$
Similar as above, $I_\la'(v_\la^1)=0$. Let $v_n^2=u_n^1-v_\la^1$, then
$$
u_n=u_\la+v_\la^1(\cdot-z_n^1)+v_n^2(\cdot-z_n^1).
$$
If $v_n^2\rg0$, i. e., $u_n^1\rg v_\la^1$ strongly in $H^1(\RN)$ as $n\rg\iy$, then
\begin{align*}
\bcs
c_\la=I_\la(u_\la)+I_\la(v_\la^1),\\
\|u_n-u_\la-v_\la^1(\cdot-z_n^1)\|\rg0,\,\hbox{as $n\rg\iy$},
\ecs
\end{align*}
and we are done. Otherwise, If $v_n^2\not\rg0$ strongly in $H^1(\RN)$ as $n\rg\iy$, similar as above,
$
\lim_{n\rg\iy}\sup_{z\in\RN}\int_{B_1(z)}|v_n^2|^2>0.
$
Then there exists $\{z_n^2\}\subset\RN$ such that $|z_n^2|\rg\iy$ as $n\rg\iy$ and
$
\lim_{n\rg\iy}\int_{B_1(z_n^2)}|v_n^2|^2>0.
$
Let $u_n^2=v_n^2(\cdot+z_n^2)$, then up to a sequence, $u_n^2\rg v_\la^2$ weakly in $H^1(\RN)$ as $n\rg\iy$ for some $v_\la^2\not=0$. Similar as above, $I_\la'(v_\la^2)=0$ and
$$
I_\la(u_n^2)\rg c_\la-I_\la(u_\la)-I_\la(v_\la^1),\,\, I_\la'(u_n^2)\rg0\,\,\hbox{in $H^{-1}(\RN)$ as $n\rg\iy$}.
$$
Let $v_n^3=u_n^2-v_\la^2$, then
$$
u_n=u_\la+v_\la^1(\cdot-z_n^1)+v_\la^2(\cdot-z_n^1-z_n^2)+v_n^3(\cdot-z_n^1-z_n^2).
$$
If $v_n^3\rg0$, i. e., $u_n^2\rg v_\la^2$ strongly in $H^1(\RN)$ as $n\rg\iy$, then
\begin{align*}
\bcs
c_\la=I_\la(u_\la)+I_\la(v_\la^1)+I_\la(v_\la^2),\\
\|u_n-u_\la-v_\la^1(\cdot-z_n^1)-v_\la^2(\cdot-z_n^1-z_n^2)\|\rg0,\,\hbox{as $n\rg\iy$},
\ecs
\end{align*}
and we are done.
Otherwise, we can repeat the procedure above. By Lemma \ref{Pohozaev}, we will have to terminate our arguments
by repeating the above proof by finite number $k$ of steps. That is, let $x_n^j=\sum_{i=1}^jz_n^i$, then
\begin{align*}
\bcs
c_\la=I_\la(u_\la)+\sum_{j=1}^kI_\la(v_\la^j),\\
\left\|u_n-u_\la-\sum_{j=1}^kv_\la^j(\cdot-x_n^j)\right\|\rg0,\,\hbox{as $n\rg\iy$}.
\ecs
\end{align*}

{\bf Step 4.} We show that after extracting a subsequence of $\{x_n^j\}$ and redefining $\{v_\la^j\}$ if necessary, Property (iii), (iv), (v) hold. Let $\lb_1, \lb_2\subset\{1,2,\cdots, k\}$ and satisfy $\lb_1\cup\lb_2=\{1,2,\cdots, k\}$ and $\{x_n^j\}_n$ is bounded if $j\in\lb_1$, $|x_n^j|\rg\iy$ as $n\rg\iy$ if $j\in\lb_2$. Then for any $j\in\lb_1$ if $\lb_1\not=\emptyset$, there exists $0\not=v^j\in H^1(\RN)$ such that, up to a sequence, $v_\la^j(\cdot-x_n^j)\rg v^j$ weakly in $H^1(\RN)$ as $n\rg\iy$ and $I_\la'(v^j)=0$ in $H^{-1}(\RN)$. By Rellich's theorem, for any $t\in[2,2N/(N-2))$, $v_\la^j(\cdot-x_n^j)\rg v^j$ strongly in $L^t(\RN)$ as $n\rg\iy$. Noting that $I_\la'(v_\la^j(\cdot-x_n^j))=0$ in $H^{-1}(\RN)$ and $I_\la(v_\la^j(\cdot-x_n^j))\le c_\la$, similar as Step 2, we know $v_\la^j(\cdot-x_n^j)\rg v^j$ strongly in $H^1(\RN)$ as $n\rg\iy$. Then, up to a sequence, there exists $\ti{v}^j\in H^1(\RN)$ such that $\sum_{j\in\lb_1}v_\la^j(\cdot-x_n^j)\rg\ti{v}^j$ strongly in $H^1(\RN)$ as $n\rg\iy$, which implies that
$
\left\|u_n-u_\la-\sum_{j\in\lb_2}v_\la^j(\cdot-x_n^j)\right\|\rg0,\,\hbox{as $n\rg\iy$}.
$
Recalling that $\|u_n-u_\la\|\not\rg0$ as $n\rg\iy$, $\lb_2\not=\emptyset$. Let $x_n^i\in\lb_2$ and
$
\lb_2^i:=\big\{j\in\lb_2: \hbox{$|x_n^i-x_n^j|$ is bounded for $n$}\big\},
$
then similar as above, up to a sequence, for some $\ti{v}_\la^i\in H^1(\RN)$, we have $\sum_{j\in\lb_2^i}v_\la^j(\cdot+x_n^i-x_n^j)\rg\ti{v}_\la^i$ strongly $H^1(\RN)$ as $n\rg\iy$. Then as $n\rg\iy$,
$
\left\|u_n-u_\la-\ti{v}_\la^i(\cdot-x_n^i)-\sum_{j\in(\lb_2\setminus\lb_2^i)}v_\la^j(\cdot-x_n^j)\right\|\rg0.
$
Without loss generality, we may assume that $\ti{v}_\la^i\not=0$. Noting that $u_n(\cdot+x_n^i)\rg\ti{v}_\la^i$ a. e. in $\RN$ as $n\rg\iy$, we get $I_\la'(\ti{v}_\la^i)=0$ in $H^{-1}(\RN)$. Then we redefine $v_\la^i:=\ti{v}_\la^i$ and as $n\rg\iy$,
$
\left\|u_n-u_\la-\sum_{j\in(\lb_2\setminus\lb_2^i)\cup\{i\}}v_\la^j(\cdot-x_n^j)\right\|\rg0.
$
By repeating the argument above by at most $(k-1)$ times and redefining $\{v_\la^j\}$ if necessary, there exists $\lb\subset\lb_2$ such that
\begin{align*}
\bcs
\hbox{$|x_n^j|\rg\iy$ and $|x_n^i-x_n^j|\rg\iy$ as $n\rg\iy$ for any $i,j\in\lb$ and $i\not=j$},\\
\|u_n-u_\la-\sum_{j\in\lb}v_\la^j(\cdot-x_n^j)\|\rg0,\,\hbox{as $n\rg\iy$}.
\ecs
\end{align*}
Finally, by Lemma \ref{Lieb}, it is easy to know $c_\la=I_\la(u_\la)+\sum_{j\in\lb}I_\la(v_\la^j)$. The proof is completed.
\qed
\vskip0.1in
\noindent{\bf Proof of Theorem \ref{Theorem 1}.} First, as a consequence of Lemma \ref{bps}, Proposition \ref{describe} and Lemma \ref{Pohozaev}, it is easy to know for almost every $\la\in J=[1/2,1]$, problem \re{q2a} admits a nontrivial solution $u_\la$ satisfying $\|u_\la\|\ge\beta$ and $\g\le I_\la(u_\la)\le c_\la$, where $\beta,\g>0$(independent of $\la$). Then there exists $\{\la_n\}\subset[1/2,1]$ and $\{u_n\}\subset H^1(\RN)$ such that as $n\rg\iy$,
\be\lab{asy}
\la_n\rg1,\,\,\g\le I_{\la_n}(u_n)\le c_{\la_n},\,\, I_{\la_n}'(u_n)=0\,\,\hbox{in $H^{-1}(\RN)$}.
\ee
By Pohoz\v{a}ev's identity,
$$
I_{\la_n}(u_n)=\frac{2+\al}{2(N+\al)}\int_{\RN}|\na u_n|^2+\frac{\al a}{2(N+\al)}\int_{\RN}|u_n|^2
$$
and $\{u_n\}$ is bounded in $H^1(\RN)$. Notice that
$$
L_a(u)=I_\la(u)+\frac{1}{2}(\la-1)\int_{\RN} (I_\al\ast F(u))F(u),\ \ u\in H^1(\RN).
$$
Then by \re{asy}, up to a sequence, there exists $c_0\in[\g,c_1]$ such that
$$
c_0:=\lim_{n\rg\iy}L_a(u_n)=\lim_{n\rg\iy}I_{\la_n}(u_n)\le\lim_{n\rg\iy}c_{\la_n}=c_1,
$$
where we used the fact that $c_\la$ is continuous from the left-hand side at $\la$. Moreover, by \re{asy}, for any $\phi\in C_0^\iy(\RN)$,
$$
\lan L_a'(u_n),\phi\ran=(\la_n-1)\int_{\RN} [I_\al\ast F(u_n)]f(u_n)\phi.
$$
Similar as above, there exists some $C>0$ such that
$$
\left(\int_{\RN}|f(u_n)\phi|^{\frac{2N}{N+\al}}\right)^{\frac{N+\al}{2N}}\le C\|\phi\|\,\,\,\mbox{uniformly for all}\,\,\phi\in C_0^\iy(\RN), n=1,2,\cdots
$$
and by the Hardy-Littlewood-Sobolev inequality,
\begin{align*}
&|\lan L_a'(u_n),\phi\ran|=(1-\la_n)\left|\int_{\RN}[I_\al\ast F(u_n)]f(u_n)\phi\right|\\
&\le C(1-\la_n)\left(\int_{\RN}|F(u_n)|^{\frac{2N}{N+\al}}\right)^{\frac{N+\al}{2N}}
\left(\int_{\RN}|f(u_n)\phi|^{\frac{2N}{N+\al}}\right)^{\frac{N+\al}{2N}}\\
&=o_n(1)\|\phi\|,
\end{align*}
where $o_n(1)\rg 0$ uniformly for any $\phi\in C_0^\iy(\RN)$ as $n\rg\iy$. That is $L_a'(u_n)\rg0$ in $H^{-1}(\RN)$ as $n\rg\iy$. In sum, we get that
$$
\|u_n\|\ge\beta,\,\,L_a(u_n)\rg c_0\le c_1,\,\, L_a'(u_n)\rg0\,\,\hbox{in $H^{-1}(\RN)$ as $n\rg\iy$}.
$$
We assume that $u_n\rg u_0$ weakly in $H^1(\RN)$ as $n\rg\iy$. If $u_n\rg u_0$ strongly in $H^1(\RN)$, then $\|u_0\|\ge\beta$ $L_a(u_0)=c_0\le c_1$ and $L_a'(u_0)=0$ in $H^{-1}(\RN)$. Otherwise, as a consequence of Proposition \ref{describe} with $\la=1, c_\la=c_0, u_\la=u_0$, there exists $k\in\mathbb{N}^+$ and $\{v^j\}_{j=1}^k\subset H^1(\RN)$ such that $v^j\not=0$, $L_a'(v^j)=0$ in $H^{-1}(\RN)$ for all $j$ and
$
c_0=L_a(u_0)+\sum_{j=1}^kL_a(v^j).
$
So let
$
\mathcal{N}:=\{u\in H^1(\RN\setminus\{0\}): L_a'(u)=0\,\,\hbox{in $H^{-1}(\RN)$}\},
$
then $\mathcal{N}\not=\emptyset$ and $\inf_{u\in\mathcal{N}}L_a(u)=E_a\in[\g,c_1]$.

\vskip0.1in
Finally, to conclude the proof of Theorem \ref{Theorem 1}, we show that $E_a$ can be achieved. Obviously, there exists $\{v_n\}\subset \mathcal{N}$ such that as $n\rg\iy$, $L_a(v_n)\rg E_a$ and $L_a'(v_n)=0$ in $H^{-1}(\RN)$. Similar as above, $\{v_n\}$ is bounded in $H^1(\RN)$. Assume that $v_n\rg v_0$ weakly in $H^1(\RN)$ as $n\rg\iy$, then $L_a'(v_0)=0$ in $H^{-1}(\RN)$. If $v_n\rg v_0$ strongly in $H^1(\RN)$, then $L_a(v_0)=E_a$. Namely, $v_0$ is a ground state solution of \re{lb1}. Otherwise, similar as above, there exists $k\in\mathbb{N}^+$ and $\{v^j\}_{j=1}^k\subset H^1(\RN)$ such that $v^j\not=0$, $L_a'(v^j)=0$ in $H^{-1}(\RN)$ for all $j$ and
$
E_a=L_a(v_0)+\sum_{j=1}^kL_a(v^j).
$
By the definition of $E_a$, $v_0=0$, $k=1$ and $L_a(v^1)=E_a$, which yields that $v^1$ is a ground state solution of \re{lb1}. The proof is completed.
\qed

\subsection{Compactness of the set of ground states solutions}

Denote the set of ground state solutions to \re{lb1} by
$$
\mathcal{N}_a:=\{u\in H^1(\RN): L_a(u)=E_a, L_a'(u)=0\,\,\hbox{in $H^{-1}(\RN)$}\},
$$
then by Theorem \ref{Theorem 1}, $\mathcal{N}_a\not=\emptyset$ for any $a>0$. Since $L_a$ is invariant by translations, $\mathcal{N}_a$ loses the compactness in $H^1(\RN)$. However, we have

\bo\lab{compact}\
For any $a>0$, up to translations, $\mathcal{N}_a$ is compact in $H^1(\RN)$.
\eo
\bp
Let $\{u_n\}\subset\mathcal{N}_a$, then $L_a(u_n)=E_a$ and $L_a'(u_n)=0$ in $H^{-1}(\RN)$. Similar as above, $\{u_n\}$ is bounded in $H^1(\RN)$. Assume that $u_n\rg u_0$ weakly in $H^1(\RN)$ as $n\rg\iy$, then $L_a'(u_0)=0$ in $H^{-1}(\RN)$. If $u_n\rg u_0$ strongly in $H^1(\RN)$, we are done. Otherwise, by virtue of Proposition \ref{describe}, up to a sequence, there exists $k\in\mathbb{N}^+$, $\{x_n^j\}_{j=1}^k\subset\RN$ and $\{v^j\}_{j=1}^k\subset H^1(\RN)$ such that $v^j\not=0$, $L_a'(v^j)=0$ in $H^{-1}(\RN)$ for all $j$ and
\begin{align*}
\bcs
E_a=L_a(u_0)+\sum_{j=1}^kL_a(v^j),\\
\|u_n-u_0-\sum_{j=1}^kv_\la^j(\cdot-x_n^j)\|\rg0\,\, \hbox{as $n\rg\iy$},
\ecs
\end{align*}
which implies that $u_0=0$, $k=1$, $v^1\in\mathcal{N}_a$ and $\|u_n(\cdot+x_n^1)-v_\la^1\|\rg0$ as $n\rg\iy$. This finishes the proof.
\ep

\subsection{Regularity, positivity and symmetry}

\noindent
Now, we adopt some ideas from \cite{MV1,alves} to give the boundedness, decay, positivity and symmetry of ground state solutions to \re{lb1}.
\bo\lab{decay}\ For any $a>0$, we have
\begin{itemize}
\item [$(i)$] $0<\inf\{\|u\|_{\infty}: u\in\mathcal{N}_a\}\le \sup\{\|u\|_{\infty}: u\in\mathcal{N}_a\}<\iy$.
\item [$(ii)$] For any $u\in \mathcal{N}_a$, $u\in C_{loc}^{1,\g}(\RN)$ for $\g\in(0,1)$.
\item [$(iii)$] For any $u\in \mathcal{N}_a$, $u$ has a constant sign and is radially symmetric about a point.
\item [$(iv)$] $E_a$ coincides with the mountain pass value.
\item [$(v)$] There exist $C,c>0$, independent of $u\in\mathcal{N}_a$, such that $|D^{\al_1}u(x)|\le C\exp(-c|x-x_0|), \,x\in \RN$ for $|\al_1|=0,1$, where $|u(x_0)|=\max_{x\in\RN}|u(x)|$.
\end{itemize}
\eo

\bp First, by Pohozaev's inequality, it is easy to know $\mathcal{N}_a$ is bounded in $H^1(\RN)$.

{\bf Claim 1.} For any $p\in[2,\frac{N}{\al}\frac{2N}{N-2})$, there exists $C_p>0$ such that
\be\lab{bounded}
\|u\|_p\le C_p\|u\|_2,\,\,\mbox{for all}\,\, u\in\mathcal{N}_a.
\ee
In fact, for any fixed $u\in\mathcal{N}_a$, let $H(u)=F(u)/u$ and $K(u)=f(u)$ in $\{x\in\RN: u(x)\not=0\}$.
Let $R>0$ and $\phi_R\in C_0^\iy(\R)$ be
such that $\phi_R(t)\in[0,1]$ for $t\in\R$, $\phi_R(t)=1$ for $|t|\le R$
and $\phi_R(t)=0$ for $|t|\ge 2R$. Setting
\begin{align*}
\bcs
H^\ast(u)=\phi_R(u)H(u),\ \ H_\ast(u)=H(u)-H^\ast(u),\\
K^\ast(u)=\phi_R(u)K(u),\ \ K_\ast(u)=K(u)-K^\ast(u).
\ecs
\end{align*}
By $(F1)$-$(F2)$, there
exists $C>0$(depending only on $R$) such that for any $x\in\RN$,
\begin{align*}
\bcs
|H^\ast(u)|\le C|u|^{\frac{\al}{N}},\ \ |K^\ast(u)|\le C|u|^{\frac{\al}{N}},\\
|H_\ast(u)|\le C|u|^{\frac{\al+2}{N-2}},\ \ |K_\ast(u)|\le C|u|^{\frac{\al+2}{N-2}}.
\ecs
\end{align*}
Obviously, $H^\ast(u), K^\ast(u)$ are uniformly bounded in $L^{2N/\al}(\RN)$ and so are $H_\ast(u), K_\ast(u)$ in $L^{2N/(\al+2)}(\RN)$ for any $u\in\mathcal{N}_a$. Thanks to the compactness of $\mathcal{N}_a$, it is easy to know for any $\e>0$, we can choose $R$ given above and depending only on $\e$ such that
$$
\left(\int_{\RN}|H_\ast(u)|^{\frac{2N}{\al+2}}\int_{\RN}|K_\ast(u)|^{\frac{2N}{\al+2}}\right)^{\frac{\al+2}{2N}}\le\e^2,\,\hbox{for all $u\in\mathcal{N}_a$}.
$$
Then repeating the argument as in \cite[Proposition 3.1]{MV1}, \re{bounded} can be concluded.
\vskip0.1in

{\bf Claim 2.} $I_\al\ast F(u)$ is uniformly bounded in $L^\iy(\RN)$ for all $u\in\mathcal{N}_a$.

\noindent By $(F1)$-$(F2)$ and the definition of the convolution of $I_\al\ast F(u)$, there exists $C(\al)$ (depending only $N,\al$) such that for any $x\in\RN$ and $u\in\mathcal{N}_a$,
{\allowdisplaybreaks
\begin{align*}
(I_\al\ast |F(u)|)(x)\le& C(\al)\int_{\R^2}(|u|^2+|u|^{(N+\al)/(N-2)})\,\ud y\\
\, &+C(\al)\int_{|x-y|\le1}\frac{|u|^2+|u|^{(N+\al)/(N-2)}}{|x-y|^{N-\al}}\,\ud y.
\end{align*}}%
Thanks to \re{bounded}, for some $c$ (independent of $u$) such that for any $x\in\RN$,
$$
(I_\al\ast |F(u)|)(x)\le c+C(\al)\int_{|x-y|\le1}\frac{|u|^2+|u|^{(N+\al)/(N-2)}}{|x-y|^{N-\al}}\,\ud y.
$$
Similar as in \cite[Proposition 2.2]{YZZ}, choosing $t\in(\frac{N}{\al},\frac{N}{\al}\frac{N}{N-2})$ with $2t\in(2,\frac{N}{\al}\frac{2N}{N-2})$ and $s\in(\frac{N}{\al},\frac{N}{\al}\frac{2N}{N+\al})$ with $s\frac{N+\al}{N-2}\in(2,\frac{N}{\al}\frac{2N}{N-2})$, there exists $C_1,C_2>0$(independent of $u$), such that
\begin{align*}
\int_{|x-y|\le1}\frac{|u|^2+|u|^{(N+\al)/(N-2)}}{|x-y|^{N-\al}}\,\ud y\le C_1\|u\|_{2t}^2+C_2\|u\|_{s\frac{N+\al}{N-2}}^{(N+\al)/(N-2)},
\end{align*}
which implies combing \re{bounded} that the claim holds.

Let $\bar{f}(x,u):=(I_\al\ast F(u))(x)f(u)$, then by $(F1)$-$(F2)$, for any $u\in\mathcal{N}_a$, $u$ satisfies that for any $\dd>0$, there exists $C_\dd>0$(independent of $u$) such that
$$
|\bar{f}(x,u)u|\le (\dd|u|^2+C_\dd|u|^{\frac{N+\al}{N-2}}),\,\, x\in\RN
$$
and
$$
-\DD u+au=\bar{f}(x,u),\,\,u\in H^1(\RN).
$$
Noting that $(N+\al)/(N-2)<2N/(N-2)$, by virtue of the standard Moser iteration \cite{GT} (see also \cite{Byeon-Zhang-Zou}), $\mathcal{N}_a$ is uniformly bounded in $L^\iy(\RN)$. Since $|\bar{f}(x,u)|=o(1)|u|$ if $\|u\|_\iy\rg0$ and $E_a>0$, it is easy to know  $\inf\{\|u\|_{\infty}: u\in\mathcal{N}_a\}>0$.
\vskip0.1in
Second, since $u\in L^\iy(\RN)$ for any $\in\mathcal{N}_a$, it follows from the elliptic estimate(see \cite{GT}) that $u\in C_{loc}^{1,\g}(\RN)$ for some $\g\in(0,1)$. From the proof of Theorem \ref{Theorem 1}, we know $E_a\le c_1$, where
$$
c_1:=\inf_{\g\in\G}\max_{t\in[0,1]}L_a(\g(t)),
$$
where $\G:=\{\g\in C([0,1], X): \g(0)=0, L_a(\g(1))<0\}$. Similar as in \cite{MV1}, for any $u\in\mathcal{N}_a$, there exists a path $\g\in\G$ such that $\g(1/2)=u$ and $L_a(\g)$ achieves its maximum at $1/2$. Thereby, $c_1=E_a$. Namely, $E_a$ is also a mountain pass value. Moreover, for any $u\in\mathcal{N}_a$, $u$ has a constant sign and is radially symmetric about some point. If $u$ is positive, then $u$ is decreasing at $r=|x-x_0|$, where $x_0$ is the maximum point of $u$. Finally, by the radial lemma, $u(x)\rg0$ uniformly as $|x-x_0|\rg\iy$ for $u\in\mathcal{N}_a$. By the comparison principle, there exist $C,c>0$, independent of $u\in\mathcal{N}_a$, such that $|D^{\al_1}u(x)|\le C\exp(-c|x-x_0|), \,x\in \RN$ for $|\al_1|=0,1$. The proof is completed.
\ep

\s{Proof of Theorem \ref{Theorem 2}}

In this section, we consider the semiclassical states of \re{q1}. To study \re{q1}, let $u(x)=v(\e x)$ and $V_{\e}(x)=V(\e x)$, then it suffices to consider the following problem
\begin{eqnarray}\label{q2}
-\DD u+V_{\e}(x)u=(I_\al\ast F(u))f(u),\,\,\,x\in\RN.
\end{eqnarray} Let $H_\e$ be the completion of $C_0^\iy(\RN)$ with respect to the norm
$$
\|u\|_\e=\left(\int_{\RN} (|\na u|^2+V_\e u^2)\right)^{\frac{1}{2}}.
$$
For any set $B\subset \RN$ and $\e>0$, we define $B_\e\equiv\{x\in\RN: \e x\in B\}$ and $B^{\dd}\equiv\{x\in\RN: \mbox{dist}(x,B)\le\dd\}$. Since we are looking for positive solutions of \re{q1}, from now on, we may assume that $f(t)=0$ for $t \le0.$  \noindent For $u\in H_\e$, let
$$
P_\e(u)=\frac{1}{2}\int_{\RN} |\na u|^2+V_\e u^2-\frac{1}{2}\int_{\RN} (I_\al\ast F(u))F(u).
$$
Fixing an arbitrary $\nu>0$, we define
\be
\chi_\e(x)=
\bcs\nonumber
0,\ \ \ \ \mbox{if}\ \ x\in O_\e,\\
\varepsilon^{-\nu},\ \ \mbox{if}\ \ x\in \RN\setminus O_\e,
\ecs
\ee
and
$$
Q_\e(u)=\left(\int_{\RN}\chi_\e u^2\, \ud x-1\right)_+^2.
$$
Let $\G_\e:H_\e\rg \R$ be given by $$\G_\e(u)=P_\e(u)+Q_\e(u).$$ To find solutions of \re{q2} which concentrate inside $O$ as $\e\rg 0$, we seek critical points $u_\e$ of $\G_\e$ satisfying $Q_\e(u_\e)=0$. The functional $Q_{\e}$ that was first introduced in \cite{BW}, will act as a penalization to force the concentration phenomena to occur inside $O$. In what follows, we seek the critical points of $\G_\e$ in some neighborhood of ground state solutions to \re{lb1} with $a=m$.

\subsection{The truncated problem}
Denote $S_m$ by the set of positive ground state solutions of \re{lb1} with $a=m$ satisfying $u(0)=\max_{x\in\RN}u(x)$, where $m$ is given in Section 1.
\bl\lab{com}
$S_m$ is compact in $H^1(\RN)$.
\el
\bp
Obviously, by Proposition \ref{decay}, $S_m\not=\emptyset$. For any $\{u_n\}\subset S_m$, without loss of generality, we assume that $u_n\rg u_0$ weakly in $H^1(\RN)$ and a. e. in $\RN$ as $n\rg\iy$. First, we claim that $u_0\not=0$. Indeed, by $(v)$ of Proposition \ref{decay}, there exist $c,C>0$ (independent of $n$) such that $|u_n(x)|\le C\exp{(-c|x|)}$ for any $x\in\RN$. By the Lebesgue dominated convergence theorem, $u_n\rg u_0$ strongly in $L^p(\RN)$ as $n\rg\iy$ for any $p\in[2,2N/(N-2)]$. So if $u_0=0$, it is easy to know $u_n\rg0$ strongly in $H^1(\RN)$ as $n\rg\iy$, which contracts the fact that $E_m>0$. Second, we claim that $u_n\rg u_0$ strongly in $H^1(\RN)$ as $n\rg\iy$. Otherwise, similar as in Proposition \ref{compact}, by Proposition \ref{describe}, there exists $k\in\mathbb{N}^+$ and $\{v^j\}_{j=1}^k\subset H^1(\RN)$ such that $v^j\not=0$, $L_m'(v^j)=0$ in $H^{-1}(\RN)$ for all $j$ and $E_m=L_m(u_0)+\sum_{j=1}^kL_m(v^j)$. Noting that $L_m(u_0)\ge E_m$ and $L_m(v^j)\ge E_m$, we get a contradiction. Finally, we show $u_0\in S_m$. Obviously, $u_0\in\mathcal{N}_m$ is positive and radially symmetric. Recalling that $0$ is the same maximum point $u_n$ for any $n$, by the local elliptic estimate, $0$ is also a maximum point of $u_0$. The proof is completed.
\ep

By Proposition \ref{decay}, let $\kappa>0$ be fixed and satisfies
\be
\lab{ee}\sup_{U\in S_m}\|U\|_\iy<\kappa.
\ee
For $k>\max_{t\in[0,\kappa]}f(t)$ fixed, let $f_k(t):=\min\{f(t),k\}$ and consider the truncated problem
\be\lab{qq8}
-\e^2\DD v+V(x)v=\e^{-\al}(I_\al\ast F_k(v))f_k(v),\ \ v\in H^1(\RN),
\ee
whose associated limit problem is
\be\lab{s2}
-\DD u+mu=(I_\al\ast F_k(u))f_k(u),\ \ u\in H^1(\RN),
\ee
where $F_k(t)=\int_0^tf_k(s)\,\ud s.$ Denote by $S_m^k$ be the set of positive ground state solutions $U$ of \eqref{s2} satisfying $U(0)=\max_{x\in \RN}U(x)$,  then by \cite[Theorem 2]{MV1}, $S_m^k\not=\emptyset$. Similar to Lemma \ref{com}, $S_m^k$ is compact in $H^1(\RN)$.
\bl\lab{l6}\
$S_m\subset S_m^k$
\el
\bp  Denote by $E_m^k$ the least energy of \re{s2}. Noting that for any $u\in S_m$, $u$ is also a solution of \re{s2}. Then we get that $E_m^k\le E_m$. By \cite{MV1}, $E_m^k$ is a mountain path value. Combing $(iv)$ of Proposition \ref{decay} and the fact that $f_k(t)\le f(t)$ for $t>0$ and $f_k(t)=f(t)=0$ for $t\le0$, we have $E_m^k\ge E_m$ and so $E_m^k=E_m$, which yields that $S_m\subset S_m^k$. The proof is completed. \ep

\subsection{Proof of Theorem \ref{Theorem 2}}

In the following, we use the truncation approach to prove Theorem \ref{Theorem 2}. Our strategy is as follows. First, we consider the truncated problem \re{qq8}. By Lemma \ref{l6}, $S_m$ is a compact subset of $S_m^k$. So we can adopt an idea in \cite{byeon} to show that \re{qq8} admits a nontrivial positive solution $v_\e$ in some neighborhood of $S_m$ for small $\e$. Second, we show that there exists $\e_0>0$ such that
$$
\|v_\e\|_\iy<\kappa,\,\hbox{for $\e\in(0,\e_0)$}.
$$
As a consequence, $v_\e$ is indeed a solution of the original problem \re{q1}.
\vskip0.1in
\noindent{\bf Completion of Proof of Theorem \ref{Theorem 2}}
\vskip0.1in
\noindent{\it Proof.}
Let
$$
\delta=\frac{1}{10}\min\{\mbox{dist}(\mathcal{M},O^c)\}.
$$
Let $\beta\in (0,\delta)$ and a cut-off $\vp\in C_0^\infty(\RN)$ such that $0\le\vp\le1,\vp(x)=1$ for $|x|\le \beta$ and $\vp(x)=0$ for $|x|\ge 2\beta$. Set $\vp_{\e}(y)=\vp(\e y), y\in\RN$ and for some $x\in (\mathcal{M})^{\beta}$ and $U\in S_m$, we define
$$
U_{\e}^x(y)=\vp_{\e}\left(y-\frac{x}{\e}\right)U\left(y-\frac{x}{\e}\right)
$$
and
$$
X_{\e}=\{U_{\e}^x\,\,|\,\,x\in (\mathcal{M})^{\beta}, U_i\in S_m\}.
$$
In the following, we show that \eqref{qq8} admits a solution in $X_\e^d$ of $X_\e$ for $\e,d>0$ small enough, where
$$
X_\e^d=\left\{u\in H_\e:  \inf_{v\in X_\e}\|u-v\|_\e\le d\right\}.
$$
In fact, since $f_k$ satisfies all the hypotheses of \cite[Theorem 2.1]{YZZ}, as a consequence, for $\e,d>0$ small, \eqref{qq8} admits a positive solution $v_{\e}\in X_\e^d$ satisfying that there exist $U\in S_m$ and a maximum point $x_\e$ of $v_\e$, such that $\lim_{\e\rg 0}dist(x_\e,\mathcal{M})=0$ and $v_\e(\e\cdot+x_\e)\rg U(\cdot+z_0)$ in $H^1(\RN)$ as $\e\rg 0$ for some $z_0\in\RN$. Noting that
$$
-\DD w_\e+V_\e(x+\frac{x_\e}{\e})w_\e=(I_\al\ast F_k(w_\e))f_k(w_\e),\ \ x\in\RN,
$$
where $w_\e(\cdot)=v_\e(\e\cdot+x_\e)$. Similar as in Proposition \ref{decay}, $I_\al\ast F_k(w_\e)$ is uniformly bounded in $L^\iy(\RN)$ for all $\e$. Then, a local elliptic estimate(see \cite{GT}) yields that $w_\e(0)\rg U(z_0)$ as $\e\rg0$. It follows from \eqref{ee} that $\|v_\e\|_\iy=w_\e(0)<\kappa$ uniformly for small $\e>0$. Therefore, for small $\e>0$, $f_k(v_\e(x))\equiv f(v_\e(x)), x\in\RN$ and then $v_\e$ is a positive solution of \re{q1}. The proof is completed.
\qed

\end{document}